\documentclass[a4paper,11pt]{amsart}  

\usepackage{dacxjo} 

\title{Lines on K3 quartic surfaces in characteristic 2}
\author{Davide Cesare Veniani}

\address{Institut für Algebraische Geometrie \\
Leibniz Universität Hannover \\
Wel\-fen\-garten 1, 30167 Hannover, Germany}
\email{veniani@math.uni-hannover.de}

\date{\today}


\begin{document}

\begin{abstract}
We prove that a K3 quartic surface defined over a field of characteristic 2 can contain at most 68 lines. If it contains 68 lines, then it is projectively equivalent to a member of a 1-dimensional family found by Rams and Sch\"utt.
\end{abstract}

\maketitle

\section{Introduction}

The enumerative geometry of lines on smooth quartic surfaces has been the object of intensive study in the last few years. 

The problem of determining the maximum number of lines on a smooth complex quartic $X$ had already been tackled in 1943 by B. Segre~\cite{segre}. However, his proof that this maximum is 64 was flawed, and it was corrected only 70 years later by Rams and Schütt~\cite{64lines}. Their proof is also valid for surfaces defined over fields $\KK$ of positive characteristic different from 2 and~3. A smooth quartic surface containing 64 lines had already been found in 1882 by Schur~\cite{schur}.

Rams and Schütt's approach is based on the study of the elliptic fibration induced by the pencil of planes containing a line on $X$. Degtyarev, Itenberg and Sertöz~\cite{AIS} -- spurred by a remark of Barth~\cite{barth} -- attacked the problem using computer-based lattice-theoretical methods. One of their main results is the complete classification of configurations of lines on smooth complex quartics with more than 52 lines. 

Quartic surfaces defined over fields of characteristic 3 and 2 feature quite different behaviors. For instance, if $\Char \KK = 3$ the Fermat surface is the unique (smooth) quartic surface up to projective equivalence which contains 112 lines. This is indeed the maximum possible, as shown by Rams and Schütt~\cite{char3}.
On the other hand, Degtyarev~\cite{degtyarev} has proven that in the case $\Char \KK = 2$ a smooth quartic surface can contain at most 60 lines, a bound attained for example by the nonsingular reduction of the surface in \cite[§8]{char2}.

The high number of lines that can be observed in characteristic~3 is due to the fact that the genus 1 fibration induced by the pencil of planes containing a line $\ell$ on $X$ can be a quasi-elliptic fibration. If this is the case, then $\ell$ intersects 30 other lines on $X$. This phenomenon, however, does not appear on smooth quartics in characteristic 2 according to \cite[Proposition 2.1]{char2} (see also Remark~\ref{rmk:smooth=>ell-lines}).

Smooth quartic surfaces are well-known examples of K3 surfaces. The minimal desingularization of a quartic surface with at most (isolated) rational double points as singularities is still a K3 surface: we call such a surface a `K3 quartic surface'. In \cite{veniani1} we extended Segre--Rams--Schütt's theorem to K3 quartic surfaces defined over fields of characteristic different from 2 and~3, providing at the same time a new proof for the smooth case which does not employ the technique of the so-called `flecnodal divisor'.

This article concerns the maximum number of lines contained on a K3 quartic surface defined over a field $\KK$ of characteristic $2$, which we can assume to be algebraically closed. The case $\Char \KK = 3$ will be dealt with in a subsequent paper. The case of complex quartic surfaces with worse singularities has been treated by González Alonso and Rams~\cite{GA-Rams}, but it is still open in positive characteristic.

Our main result is the following theorem.

\begin{theorem} \label{thm:char2}
    If $X$ is a K3 quartic surface defined over a field of characteristic~$2$, then $X$ contains at most 68 lines.
\end{theorem}

The striking fact that non-smooth K3 quartic surfaces can contain more lines than smooth surfaces was already noted by Rams and Schütt~\cite{char2}, who also provided a 1-dimensional family $\XXX$ of quartic surfaces with one singularity of type~$\bA_3$ attaining the bound of Theorem~\ref{thm:char2}.
This feature can be explained by the fact that the phenomenon of quasi-elliptic lines (i.e., lines inducing a quasi-elliptic fibration on the minimal desingularization of $X$) does indeed take place once we allow for rational double points.

We are also able to prove the following uniqueness result.
\begin{theorem} \label{thm:char2-uniqueness}
    If $X$ contains 68 lines, then $X$ is projectively equivalent to a member of family $\XXX$.
\end{theorem}

The paper is structured as follows.
\begin{description}
	\item[Section~\ref{sec:K3-quartic}] We set the notation and present some general facts on K3 quartic surfaces defined over a field of any characteristic.
    \item[Section~\ref{sec:char2-elliptic}] We list the main results about elliptic lines, whose study is essentially the same as in characteristic~0.
    \item[Section~\ref{sec:char2-qe}] Quasi-elliptic lines, which provide the main differences both from the case of fields of characteristic different from 2, and from the case of smooth surfaces, are studied extensively in this section.
    \item[Section~\ref{sec:char2-proof}] We carry out the proof of Theorem~\ref{thm:char2}.
    \item[Section~\ref{sec:char2-family}] We carry out the proof of Theorem~\ref{thm:char2-uniqueness}.
\end{description}

\subsection*{Acknowledgments}

I wish to thank my supervisor Matthias Schütt and Víctor González Alonso for helpful comments and discussions.

\section{K3 quartic surfaces} \label{sec:K3-quartic}

We collect here some general facts about K3 quartic surfaces. Most of them do not depend on the characteristic of the ground field, so in this section we work over a fixed algebraically closed ground field $\KK$ of characteristic $p \geq 0$.

\begin{definition} \label{def:K3quartic}
 A \emph{K3 quartic surface} is a surface in $\PP^3$ of degree~$4$ with at most rational double points as singularities.
\end{definition}

Henceforth, $X$ will denote a fixed K3 quartic surface, $\Sing(X)$ the set of singular points on $X$ (which are isolated because rational double points are normal singularities) and $\rho : Z \rightarrow X$ the minimal desingularization of $X$.

The name `K3 quartic surface' stems from the fact that $Z$ is a K3 surface. By standard arguments in the theory of surfaces, the line bundle $H := \rho^*(\mathcal O_{X}(1))$ is a nef line bundle such that $H^2 = 4$.

\subsection{Lines and singular points}

From now on, $\ell$ will denote a line lying on a K3 quartic surface $X$ with minimal desingularization $\rho: Z \rightarrow X$. Any divisor in the complete linear system defined by $H$ will be called a \emph{hyperplane divisor} (and often denoted by $H$, too). The strict transform of $\ell$ will be denoted by $L$.

\begin{lemma} \label{lem:genus-1-fib}
 The pencil of planes $\{\Pi_t\}_{t \in \PP^1}$ containing the line $\ell$ induces a genus~1 fibration $\pi: Z \rightarrow \PP^1$.
\end{lemma}
\proof
The proof is similar to \cite[Lemma 1.2]{veniani1}, although the generic fiber is not necessarily smooth if $\Char \KK =2$ or $3$.
\endproof

The genus 1 fibration $\pi$ induced by the pencil of planes containing $\ell$ -- or, for short, induced by $\ell$ -- is always elliptic if $\Char \KK \neq 2,3$, but in general we have to distinguish the two cases. 

\begin{definition}
A line $\ell$ is said to be \emph{elliptic} (respectively \emph{quasi-elliptic}) if it induces an elliptic (respectively \emph{quasi-elliptic}) fibration.
\end{definition}

A fiber $F_t$ of $\pi$ ($t \in \PP^1$) is the pullback through $\rho$ of the residual cubic $E_t$ contained in $\Pi_t$. 
We denote the restriction of $\pi$ to $L$ again by $\pi$.

\begin{definition}
    If the morphism $\pi: L \rightarrow \PP^1$ is constant, we say that $L$ has degree $0$; otherwise, the \emph{degree} of $\ell$ is the degree of the morphism $\pi: L \rightarrow \PP^1$.
\end{definition}

\begin{definition}
    The \emph{singularity} of a line $\ell$ is the number of singular points of $X$ lying on $\ell$.
\end{definition}

As shown in~\cite[Proposition 1.5]{veniani1}, the degree of a line is never greater than~$3$; therefore, the morphism $\pi:L\rightarrow \PP^1$ is always separable if $\Char \KK \neq 2,3$. Again, in general we have to make a distinction.

\begin{definition}
    A line $\ell$ is said to be \emph{separable} (respectively \emph{inseparable}) if the induced morphism $\pi:L\rightarrow \PP^1$ is separable (respectively inseparable).
\end{definition}

Given a separable line $\ell$, we will say that a point $P$ on $\ell$ is a point of ramification $n_m$ if the corresponding point on $L$ has ramification index $n$ and $\length(\Omega_{L/\PP^1}) = m$. We recall that if $\Char \KK$ does not divide $n$, then $m = n-1$ and can be omitted, whereas if $\Char \KK$ divides $n$, then $m \geq n$. 

\begin{lemma} \label{lemma:linesthroughsingularpoint}
    If $P$ is a singular point on a K3 quartic surface $X$, then there are at most $8$ lines lying on $X$ and passing through $P$.
\end{lemma}
\proof
See~\cite[Lemma 1.7]{veniani1}.
\endproof

\begin{lemma} \label{lemma:gen-res-cub-smooth}
    If $P$ is a singular point on a line $\ell$, a general residual cubic relative to $\ell$ is smooth at $P$.
\end{lemma}
\proof
If this is not the case, then $P$ is a triple point, as can be checked by an explicit computation. \endproof

\begin{lemma} \label{lemma:exc-div-sec-or-fib-comp}
    Let $P$ be a singular point on a line $\ell$. Then, exactly one of the exceptional divisors on $Z$ coming from $P$ is a section of the fibration induced by $\ell$, and all others are fiber components.
\end{lemma}
\proof
Take a general residual cubic $E$ relative to $\ell$. Since $E$ is smooth at $P$ by Lemma \ref{lemma:gen-res-cub-smooth}, its strict transform $F$ hits exactly one exceptional divisor. On the other hand, all other exceptional divisors have intersection $0$ with the fiber $F$, so they must be fiber components.
\endproof

\subsection{Valency}

\begin{definition}
    Given a K3 quartic surface $X$, we will denote by $\Phi(X)$ the number of lines lying on $X$.
\end{definition}

Building on an idea of Segre~\cite{segre}, we will usually be interested in finding a \emph{completely reducible plane}, i.e., a plane $\Pi$ such that the intersection $X \cap \Pi$ splits into the highest possible number of irreducible components, namely four lines $\ell_1,\ldots,\ell_4$ (not necessarily distinct). If a line $\ell'$ not lying on $\Pi$ meets two or more distinct lines $\ell_i$, then their point of intersection must be a singular point of the surface. It follows that $\Phi(X)$ is bounded by
\begin{equation} \label{eq:PhiX}
\begin{split}
    \Phi(X) &\leq \#\{\text{lines in $\Pi$}\} \\
            &\quad + \#\{\text{lines not in $\Pi$ going through $\Pi\cap\Sing(X)$} \} \\ 
            &\quad + \sum_{i = 1}^4\#\{\text{lines not in $\Pi$ meeting $\ell_i$ in a smooth point}\}. \\
\end{split}
\end{equation}

It will then be a matter of finding a bound for the second and third contribution. The former will be usually dealt with using Lemma \ref{lemma:linesthroughsingularpoint}. As for the latter, it is natural to introduce the following definition.
\begin{definition}
    The \emph{valency} of $\ell$, denoted by $v(\ell)$, is the number of lines on $X$ distinct from $\ell$ which intersect $\ell$ in smooth points.
\end{definition}

Most of the time we will express the latter contribution in terms of $v(\ell_i)$, and much of the work will be dedicated to finding a bound for these quantities. 
Of course, not all K3 quartic surfaces admit a completely reducible plane, in which case we will turn to other techniques, such as the ones presented in \cite[§5]{veniani1}.

\begin{definition}
    A \emph{3-fiber} is a fiber whose residual cubic splits into three lines, whereas a \emph{1-fiber} is a fiber whose residual cubic splits into a line and an irreducible conic. 
    A line $\ell$ is said to be \emph{of type $(p,q)$}, $p,q\geq 0$, if in its fibration there are $p$ fibers of the former kind and $q$ fibers of the latter kind.
\end{definition}

\begin{definition} The (local) valency of a fiber $F$, denoted by $v_\ell(F)$, is the number of lines distinct from $\ell$ contained in the plane corresponding to $F$ that meet $\ell$ in a smooth point. When it is clear from the context, we will simply write $v(F)$.
\end{definition}

Obviously, 
\begin{equation} \label{eq:v=sum-vF}
v(\ell) = \sum_{t \in \PP^1} v_\ell(F_t),
\end{equation}
and the sum is actually a finite sum. 

\begin{lemma} \label{lemma:no-p-fibers}
    If $\ell$ is an elliptic line without 3-fibers, then $v(\ell) \leq 12$.
\end{lemma}
\proof See~\cite[Lemma 2.7]{veniani1}. \endproof

The assumption that $\ell$ is elliptic is essential, see §\ref{sec:char2-qe}.

\subsection{Lines of the first and second kind}

Let $x_0,x_1,x_2,x_3$ be the coordinates of $\PP^3$. Up to projective equivalence, we can suppose that the line $\ell$ is given by the vanishing of $x_0$ and $x_1$, so that the quartic $X$ is defined by
\begin{equation} \label{eq:surfaceX}
 X: \sum_{i_0 + i_1 + i_2 + i_3 = 4} a_{i_0 i_1 i_2 i_3} x_{0}^{i_0} x_{1}^{i_1} x_{2}^{i_2} x_{3}^{i_3} = 0, \quad a_{i_0 i_1 0 0} = 0 \text{ for all $i_0$, $i_1$},
\end{equation}
where $i_0,\ldots,i_4$ are non-negative integers. 

\begin{remark} \label{rmk:res-cub}
We will usually parametrize the planes containing $\ell$ by $\Pi_t: x_0 = t x_1$, $t\in \PP^1$, where of course $t = \infty$ denotes the plane $x_1 = 0$. Two equations which define the \emph{residual cubic} $E_t$ contained in $\Pi_t$ are the equation of $\Pi_t$ itself and the equation $g \in \KK[t][x_1,x_2,x_3]_{(3)}$ obtained by substituting $x_0$ with $t x_1$ in \eqref{eq:surfaceX} and factoring out $x_1$. An explicit computation shows that the intersection of $\ell$ with $E_t$ is given by the points $[0:0:x_2:x_3]$ satisfying
\begin{equation} \label{eq:penciloncurve}
 g_t(0,x_2,x_3) = t \alpha(x_2, x_3) + \beta(x_2, x_3) = 0.
\end{equation}
\end{remark}

Given a line $\ell$ of positive degree, a crucial technique to find bounds for $v(\ell)$ is to count the points of intersection of the residual cubics $E_t$ and $\ell$ which are inflection points for $E_t$. By \emph{inflection point} we mean here a point which is also a zero of the hessian of the cubic. In fact, if a residual cubic $E_t$ contains a line as a component, all the points of the line will be inflection points of $E_t$. 

Supposing that the surface $X$ is defined as in equation \eqref{eq:surfaceX}, the hessian of the equation $g$ defining the residual cubic $E_t$ (see Remark \ref{rmk:res-cub}) restricted on the line $\ell$ is given by
\begin{equation} \label{eq:hessian}
h := \det \left( \frac{\partial^2 g}{\partial x_i x_j} \right)_{1\leq i,j, \leq 3} \bigg\rvert_{x_1 = 0} \in \KK[t][x_2,x_3]_{(3)},
\end{equation}
which is a polynomial of degree $5$ in $t$, with forms of degree $3$ in $(x_2,x_3)$ as coefficients. 
If $\Char \KK = 2$, we need to modify the definition of the hessian slightly, as suggested by Rams and Schütt \cite{char3}. If $m$ is the coefficient of the monomial $x_1 x_2 x_3$ in $g$, then one defines
\begin{equation} \label{eq:hessian-char2}
    \tilde h = \frac14 \bigg( \frac 12 \det \Big( \frac{\partial^2 g}{\partial x_i x_j} \Big)_{1\leq i,j, \leq 3} - m^2 g \bigg) \bigg\rvert_{x_1 = 0} \in \KK[t][x_2,x_3]_{(3)},
\end{equation}
which is to be understood first as an algebraic expression over $\ZZ$ in terms of the generic coefficients of $g$, then interpreted over $\KK$ by reducing modulo 2 and substituting.

We want now to find the number of lines intersecting $\ell$ by studying the common solutions of \eqref{eq:penciloncurve} and \eqref{eq:hessian} (or \eqref{eq:hessian-char2}) on the line $\ell$. It is convenient to extend Segre's nomenclature~\cite{segre}.

\begin{definition}
    The resultant $R(\ell)$ with respect to the variable $t$ of the polynomials \eqref{eq:penciloncurve} and \eqref{eq:hessian} (or \eqref{eq:hessian-char2} if $\Char \KK = 2$) is called the \emph{resultant} of the line $\ell$.
\end{definition}

\begin{definition}
 We say that a line $\ell$ of positive degree is a line of the \emph{second kind} if its resultant is identically equal to zero. Otherwise, we say that $\ell$ is a line of the \emph{first kind}.
\end{definition}

A root $[\bar x_2: \bar x_3]$ of $R(\ell)$ corresponds to a point $P = [0:0:\bar x_2: \bar x_3]$ on $\ell$; if $P$ is a smooth surface point, then it is an inflection point for the residual cubic passing through it.
A local computation yields the following lemma, which holds in any characteristic.

\begin{lemma} \label{lem:1stkind-loc}
    Let $\ell$ be a line of the first kind. 
    \begin{enumerate}[(a)]
        \item If one line intersects $\ell$ at a smooth point $P$, then $P$ is a root of $R(\ell)$.
        \item If two lines or one double line intersect $\ell$ at a smooth point $P$, then $P$ is a root of $R(\ell)$ of order at least $2$.
        \item If three lines, one double line and a simple line, or one triple line intersect $\ell$ at a smooth point $P$, then $P$ is a root of $R(\ell)$ of order at least $5$.
    \end{enumerate}
\end{lemma}

\begin{proposition} \label{prop:1stkind}
If $\ell$ is a line of the first kind, then $v(\ell) \leq 3 + 5\,d$.
\end{proposition}
\proof Since equation \eqref{eq:penciloncurve} is linear in $t$ and -- once one has got rid of the common factors of $\alpha$ and $\beta$ -- it has degree $d$ in $(x_2,x_3)$, the resultant $R(\ell)$ of a line $\ell$ of the first kind is a form in $(x_2,x_3)$ of degree $3+5\,d$. The claim follows from Lemma \ref{lem:1stkind-loc}.
\endproof

\subsection{Triangle free surfaces} \label{subsec:triangle-free}

We follow here the nomenclature of \cite[§5]{veniani1}. In particular, the \emph{line graph} of a K3 quartic surface $X$ is the dual graph of the strict transforms of its lines on $Z$. The line graph $\Gamma = \Gamma(X)$ of a K3 quartic surface $X$ is a graph without loops or multiple edges. By definition, the number of its vertices is equal to $\Phi(X)$.

A Dynkin diagram (resp. extended Dynkin diagram) is also called an \emph{elliptic} graph (resp. \emph{parabolic} graph).

\begin{definition}
    A K3 quartic surface $X$ is called \emph{triangle free} if its line graph contains no triangles, i.e., cycles of length~$3$.
\end{definition}

In other words, a K3 quartic surface $X$ is triangle free if there are no triples of lines on $X$ intersecting pairwise in \emph{smooth} points. The next definition has an analogous geometric interpretation.

\begin{definition}
    A K3 quartic surface $X$ is called \emph{square free} if it is triangle free and if its line graph contains no squares, i.e., cycles of length~$4$.
\end{definition}

\begin{lemma} \label{lem:elliptic-triangle-free}
    If $\ell$ is an elliptic line on a triangle free K3 quartic surface, then $v(\ell) \leq 12$.
\end{lemma}
\proof
For lines of degree 0, we have $v(\ell) \leq 2$ by \cite[Lemma 2.6]{veniani1}, so we can suppose that $\ell$ has positive degree. It is sufficient to show that 
\begin{equation} \label{eq:triangle-free}
    2\, v(F) \leq e(F) 
\end{equation}
for all fibers $F$ of the fibration induced by $\ell$.
In fact, if this holds, then by formula \eqref{eq:v=sum-vF}
\[
    2\,v(\ell) = \sum_{t \in \PP^1} 2\,v(F_t) \leq \sum_{t \in \PP^1} e(F_t) \leq 24
\]
and we conclude.

Formula \eqref{eq:triangle-free} is clear for a 1-fiber $F$, because $v(F) \leq 1$ and $e(F) \geq 2$.
Observe that if a 3-fiber $F$ contains a double line, then $e(F) \geq 6$. Since $v(F)$ can never be greater than $3$, we get \eqref{eq:triangle-free}. Hence, we can suppose that $F$ is a 3-fiber composed of three distinct lines $\ell_1$, $\ell_2$, $\ell_3$. Moreover, no three lines among $\ell$ and the $\ell_i$'s can meet in a smooth point, since there are no triangles. 

If $v(F) = 3$, then two configurations may arise, as pictured below: either the $\ell_i$'s meet in different points or they are concurrent. If they meet in different points, then all points must be singular, because of the triangle free assumption, giving rise to a fiber of type $\I_n$ with $n\geq 6$. If they are concurrent, then the intersection point must be singular and the corresponding fiber must have at least 4 components, three of which of multiplicity 1, and no cycle, i.e., it must be of type $\I_n^*$ or $\IV^*$. In any case, $e(F_s) \geq 6$ and again we obtain \eqref{eq:triangle-free}.

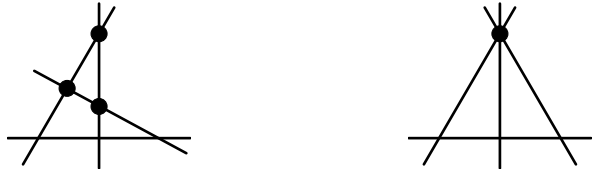
\begin{figure}[ht]
\centering
\begin{tikzpicture}[line cap=round,line join=round,>=triangle 45,x=0.4cm,y=0.4cm]
    \begin{scope}[xshift=0,yshift=0]
        \clip(-3,-1) rectangle (3,5);
        \draw [line width=1pt] (-2.5,-0.87)-- (0.5,4.33);
        \draw [line width=1pt] (-2.13,2.23)-- (2.87,-0.5);
        \draw [line width=1pt] (-3,0)-- (3,0);
        \draw [line width=1pt] (0,4.46)-- (0,-1);
        \filldraw [color=black] (0,3.46) circle (3pt);
        \filldraw [color=black] (0,1.05) circle (3pt);
        \filldraw [color=black] (-1.05,1.65) circle (3pt);
    \end{scope}

    \begin{scope}[xshift=150,yshift=0]
        \clip(-3,-1) rectangle (3,5);
        \draw [line width=1pt] (-2.5,-0.87)-- (0.5,4.33);
        \draw [line width=1pt] (2.5,-0.87)-- (-0.5,4.33);
        \draw [line width=1pt] (-3,0)-- (3,0);
        \draw [line width=1pt] (0,4.46)-- (0,-1);
        \filldraw [color=black] (0,3.46) circle (3pt);
    \end{scope}
\end{tikzpicture}
\end{figure}

If $v(F) = 2$, then there is a singular point on $\ell$ and one of the $\ell_i$'s passes through it, while the other two lines meet $\ell$ in two other smooth points. Again, the lines $\ell_i$ can meet in different points or in the same point (see picture below). In the former case, it is not possible that all the intersection points of the $\ell_i$'s are smooth; thus, the fiber is of type $\I_n$ with $n\geq 4$ and \eqref{eq:triangle-free} holds. In the latter case, we can argue as before.

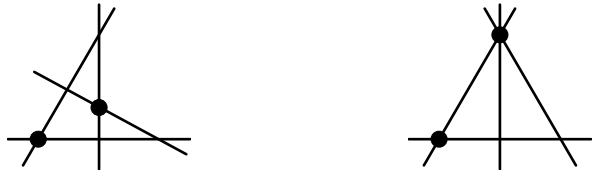
\begin{figure}[ht]
\centering
\begin{tikzpicture}[line cap=round,line join=round,>=triangle 45,x=0.4cm,y=0.4cm]
    \begin{scope}[xshift=0,yshift=0]
        \clip(-3,-1) rectangle (3,5);
        \draw [line width=1pt] (-2.5,-0.87)-- (0.5,4.33);
        \draw [line width=1pt] (-2.13,2.23)-- (2.87,-0.5);
        \draw [line width=1pt] (-3,0)-- (3,0);
        \draw [line width=1pt] (0,4.46)-- (0,-1);
        \filldraw [color=black] (-2,0) circle (3pt);
        \filldraw [color=black] (0,1.05) circle (3pt);
    \end{scope}

    \begin{scope}[xshift=150,yshift=0]
        \clip(-3,-1) rectangle (3,5);
        \draw [line width=1pt] (-2.5,-0.87)-- (0.5,4.33);
        \draw [line width=1pt] (2.5,-0.87)-- (-0.5,4.33);
        \draw [line width=1pt] (-3,0)-- (3,0);
        \draw [line width=1pt] (0,4.46)-- (0,-1);
        \filldraw [color=black] (-2,0) circle (3pt);
        \filldraw [color=black] (0,3.46) circle (3pt);
    \end{scope}
\end{tikzpicture}
\end{figure}

In case $v(F) = 1$ formula \eqref{eq:triangle-free} is automatically satisfied, since for a 3-fiber $e(F) \geq 3$.
\endproof

In the last part of this section, we would like to classify the possible configurations of lines and singular points on a completely reducible plane. Note that if three lines form a triangle, then they are necessarily coplanar and the plane containing them is completely reducible.

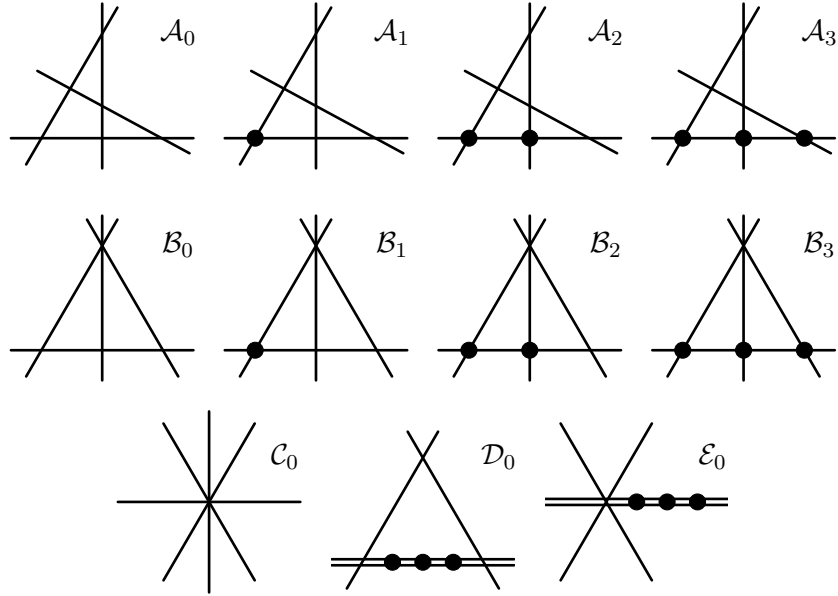
\begin{figure}
\vspace{0.5cm}
\centering
\begin{tikzpicture}[line cap=round,line join=round,>=triangle 45,x=0.4cm,y=0.4cm]
    \begin{scope}[xshift=0,yshift=0]
        \draw [line width=1pt] (-2.5,-0.87)-- (0.5,4.33);
        \draw [line width=1pt] (-2.13,2.23)-- (2.87,-0.5);
        \draw [line width=1pt] (-3,0)-- (3,0);
        \draw [line width=1pt] (0,4.46)-- (0,-1);
        \draw (2.5,3.5) node {$\mathcal A_0$};
    \end{scope}
    \begin{scope}[xshift=80,yshift=0]
        \draw [line width=1pt] (-2.5,-0.87)-- (0.5,4.33);
        \draw [line width=1pt] (-2.13,2.23)-- (2.87,-0.5);
        \draw [line width=1pt] (-3,0)-- (3,0);
        \draw [line width=1pt] (0,4.46)-- (0,-1);
        \filldraw [color=black] (-2,0) circle (3pt);
        \draw (2.5,3.5) node {$\mathcal A_1$};
    \end{scope}
    \begin{scope}[xshift=160,yshift=0]
        \draw [line width=1pt] (-2.5,-0.87)-- (0.5,4.33);
        \draw [line width=1pt] (-2.13,2.23)-- (2.87,-0.5);
        \draw [line width=1pt] (-3,0)-- (3,0);
        \draw [line width=1pt] (0,4.46)-- (0,-1);
        \filldraw [color=black] (-2,0) circle (3pt);
        \filldraw [color=black] (0,0) circle (3pt);
        \draw (2.5,3.5) node {$\mathcal A_2$};
    \end{scope}
    \begin{scope}[xshift=240,yshift=0]
        \draw [line width=1pt] (-2.5,-0.87)-- (0.5,4.33);
        \draw [line width=1pt] (-2.13,2.23)-- (2.87,-0.5);
        \draw [line width=1pt] (-3,0)-- (3,0);
        \draw [line width=1pt] (0,4.46)-- (0,-1);
        \filldraw [color=black] (-2,0) circle (3pt);
        \filldraw [color=black] (0,0) circle (3pt);
        \filldraw [color=black] (2,0) circle (3pt);
        \draw[color=black] (2.5,3.5) node {$\mathcal A_3$};
    \end{scope}

    \begin{scope}[xshift=0,yshift=-80]
        \draw [line width=1pt] (-2.5,-0.87)-- (0.5,4.33);
        \draw [line width=1pt] (2.5,-0.87)-- (-0.5,4.33);
        \draw [line width=1pt] (-3,0)-- (3,0);
        \draw [line width=1pt] (0,4.46)-- (0,-1);
        \draw[color=black] (2.5,3.5) node {$\mathcal B_0$};
    \end{scope}
    \begin{scope}[xshift=80,yshift=-80]
        \draw [line width=1pt] (-2.5,-0.87)-- (0.5,4.33);
        \draw [line width=1pt] (2.5,-0.87)-- (-0.5,4.33);
        \draw [line width=1pt] (-3,0)-- (3,0);
        \draw [line width=1pt] (0,4.46)-- (0,-1);
        \filldraw [color=black] (-2,0) circle (3pt);
        \draw (2.5,3.5) node {$\mathcal B_1$};
    \end{scope}
    \begin{scope}[xshift=160,yshift=-80]
        \draw [line width=1pt] (-2.5,-0.87)-- (0.5,4.33);
        \draw [line width=1pt] (2.5,-0.87)-- (-0.5,4.33);
        \draw [line width=1pt] (-3,0)-- (3,0);
        \draw [line width=1pt] (0,4.46)-- (0,-1);
        \filldraw (-2,0) circle (3pt);
        \filldraw (0,0) circle (3pt);
        \draw (2.5,3.5) node {$\mathcal B_2$};
    \end{scope}
    \begin{scope}[xshift=240,yshift=-80]
        \draw [line width=1pt] (-2.5,-0.87)-- (0.5,4.33);
        \draw [line width=1pt] (2.5,-0.87)-- (-0.5,4.33);
        \draw [line width=1pt] (-3,0)-- (3,0);
        \draw [line width=1pt] (0,4.46)-- (0,-1);
        \filldraw (-2,0) circle (3pt);
        \filldraw (0,0) circle (3pt);
        \filldraw (2,0) circle (3pt);
        \draw[color=black] (2.5,3.5) node {$\mathcal B_3$};
    \end{scope}
    
    \begin{scope}[xshift=40,yshift=-160]
        \draw [line width=1pt] (-1.5, -0.6)-- (1.5, 4.6);
        \draw [line width=1pt] (-1.5,4.6)-- (1.5,-0.6);
        \draw [line width=1pt] (-3,2)-- (3,2);
        \draw [line width=1pt] (0,5)-- (0,-1);
        \draw (2.5,3.5) node {$\mathcal C_0$};
    \end{scope}
    
        \begin{scope}[xshift=120,yshift=-160]
        \clip(-3,-1) rectangle (3,5);
        \draw [line width=1pt] (-2.5,-0.87)-- (0.5,4.33);
        \draw [line width=1pt] (2.5,-0.87)-- (-0.5,4.33);
        \draw [line width=1pt] (-3,-0.1)-- (3,-0.1);
        \draw [line width=1pt] (-3,0.1)--(3,0.1);
        \draw[color=black] (2.5,3.5) node {$\mathcal D_0$};
        \filldraw [color=black] (-1,0) circle (3pt);
        \filldraw [color=black] (0,0) circle (3pt);
        \filldraw [color=black] (1,0) circle (3pt);
    \end{scope}
    \begin{scope}[xshift=200,yshift=-160]
        \clip(-3,-1) rectangle (3,5);
        \draw [line width=1pt] (-2.5, -0.6)-- (0.5, 4.6);
        \draw [line width=1pt] (-2.5,4.6)-- (0.5,-0.6);
        \draw [line width=1pt] (-3,1.9)-- (3,1.9);
        \draw [line width=1pt] (-3,2.1)-- (3,2.1);
        \draw[color=black] (2.5,3.5) node {$\mathcal E_0$};
        \filldraw [color=black] (0,2) circle (3pt);
        \filldraw [color=black] (1,2) circle (3pt);
        \filldraw [color=black] (2,2) circle (3pt);
    \end{scope}
\end{tikzpicture}
\caption{\small Possible configurations of lines on a plane with a triangle. Singular points are marked with a bullet. In configurations $\mathcal D_0$ and $\mathcal E_0$ the singular points might coincide.}
\label{fig:triangle-conf}
\end{figure}

\begin{lemma} \label{lem:triangle-conf}
If three lines on $X$ form a triangle, then they are contained in plane $\Pi$ such that the intersection of $\Pi$ and $X$ has one of the configurations pictured in Figure~\ref{fig:triangle-conf}.
\end{lemma}
\proof
Let $\ell_1, \ell_2, \ell_3$ be the lines forming a triangle and $\ell_4$ the fourth line on the plane. If $\ell_4$ coincides with one of the former, then we get configurations $\mathcal D_0$ or $\mathcal E_0$. Suppose the four lines are pairwise distinct. A priori the following three configurations are possible:
\begin{itemize}
\item either the lines meet in pairwise distinct point (configurations $\mathcal A$),
\item or exactly three of them are concurrent (configurations $\mathcal B$),
\item or four of them are concurrent (configuration $\mathcal C$).
\end{itemize}

Note that a singular point of the surface contained in the plane must be the intersection point of two or more lines. By hypothesis, $\ell_1$, $\ell_2$ and $\ell_3$ meet in smooth points. Up to symmetry, the only possible configurations are those in the picture.
\endproof

\begin{figure}
\vspace{0.5cm}
\centering
\begin{tikzpicture}[line cap=round,line join=round,>=triangle 45,x=0.4cm,y=0.4cm]
    \begin{scope}[xshift=0,yshift=0]
        \draw [line width=1pt] (-2.5,-0.87)-- (0.5,4.33);
        \draw [line width=1pt] (-2.13,2.23)-- (2.87,-0.5);
        \draw [line width=1pt] (-3,0)-- (3,0);
        \draw [line width=1pt] (0,4.46)-- (0,-1);
        \filldraw (0,3.46) circle (3pt);
        \filldraw (0,1.05) circle (3pt);
        \filldraw (-1.05,1.65) circle (3pt);
        \draw (2.5,3.5) node {$\mathcal A_4$};
    \end{scope}
    \begin{scope}[xshift=80,yshift=0]
        \draw [line width=1pt] (-2.5,-0.87)-- (0.5,4.33);
        \draw [line width=1pt] (-2.13,2.23)-- (2.87,-0.5);
        \draw [line width=1pt] (-3,0)-- (3,0);
        \draw [line width=1pt] (0,4.46)-- (0,-1);
        \filldraw  (-2,0) circle (3pt);
        \filldraw (0,1.05) circle (3pt);
        \draw (2.5,3.5) node {$\mathcal A_5$};
    \end{scope}
    \begin{scope}[xshift=160,yshift=0]
        \draw [line width=1pt] (-2.5,-0.87)-- (0.5,4.33);
        \draw [line width=1pt] (-2.13,2.23)-- (2.87,-0.5);
        \draw [line width=1pt] (-3,0)-- (3,0);
        \draw [line width=1pt] (0,4.46)-- (0,-1);
        \filldraw  (-2,0) circle (3pt);
        \filldraw (0,1.05) circle (3pt);
        \filldraw (-1.05,1.65) circle (3pt);
        \draw (2.5,3.5) node {$\mathcal A_6$};
    \end{scope}
    \begin{scope}[xshift=240,yshift=0]
        \draw [line width=1pt] (-2.5,-0.87)-- (0.5,4.33);
        \draw [line width=1pt] (-2.13,2.23)-- (2.87,-0.5);
        \draw [line width=1pt] (-3,0)-- (3,0);
        \draw [line width=1pt] (0,4.46)-- (0,-1);
        \filldraw  (-2,0) circle (3pt);
        \filldraw (0,3.46) circle (3pt);
        \filldraw (0,1.05) circle (3pt);
        \filldraw (-1.05,1.65) circle (3pt);
        \draw (2.5,3.5) node {$\mathcal A_7$};
    \end{scope}
    \begin{scope}[xshift=0,yshift=-80]
        \draw [line width=1pt] (-2.5,-0.87)-- (0.5,4.33);
        \draw [line width=1pt] (-2.13,2.23)-- (2.87,-0.5);
        \draw [line width=1pt] (-3,0)-- (3,0);
        \draw [line width=1pt] (0,4.46)-- (0,-1);
        \filldraw  (-2,0) circle (3pt);
        \filldraw  (0,0) circle (3pt);
        \filldraw (0,1.05) circle (3pt);
        \filldraw (-1.05,1.65) circle (3pt);
        \draw (2.5,3.5) node {$\mathcal A_8$};
    \end{scope}
    \begin{scope}[xshift=80,yshift=-80]
        \draw [line width=1pt] (-2.5,-0.87)-- (0.5,4.33);
        \draw [line width=1pt] (-2.13,2.23)-- (2.87,-0.5);
        \draw [line width=1pt] (-3,0)-- (3,0);
        \draw [line width=1pt] (0,4.46)-- (0,-1);
        \filldraw  (-2,0) circle (3pt);
        \filldraw  (0,0) circle (3pt);
        \filldraw (0,3.46) circle (3pt);
        \filldraw (0,1.05) circle (3pt);
        \filldraw (-1.05,1.65) circle (3pt);
        \draw (2.5,3.5) node {$\mathcal A_9$};
    \end{scope}
    \begin{scope}[xshift=160,yshift=-80]
        \draw [line width=1pt] (-2.5,-0.87)-- (0.5,4.33);
        \draw [line width=1pt] (-2.13,2.23)-- (2.87,-0.5);
        \draw [line width=1pt] (-3,0)-- (3,0);
        \draw [line width=1pt] (0,4.46)-- (0,-1);
        \filldraw  (-2,0) circle (3pt);
        \filldraw  (0,0) circle (3pt);
        \filldraw  (2,0) circle (3pt);
        \filldraw (0,3.46) circle (3pt);
        \filldraw (0,1.05) circle (3pt);
        \filldraw (-1.05,1.65) circle (3pt);
        \draw (2.5,3.5) node {$\mathcal A_{10}$};
    \end{scope}

    \begin{scope}[xshift=240,yshift=-80]
        \draw [line width=1pt] (-2.5,-0.87)-- (0.5,4.33);
        \draw [line width=1pt] (2.5,-0.87)-- (-0.5,4.33);
        \draw [line width=1pt] (-3,0)-- (3,0);
        \draw [line width=1pt] (0,4.46)-- (0,-1);
        \filldraw (0,3.46) circle (3pt);
        \draw (2.5,3.5) node {$\mathcal B_4$};
    \end{scope}
    
    \begin{scope}[xshift=0,yshift=-160]
        \clip(-3,-1) rectangle (3,5);
        \draw [line width=1pt] (-2.5,-0.87)-- (0.5,4.33);
        \draw [line width=1pt] (2.5,-0.87)-- (-0.5,4.33);
        \draw [line width=1pt] (-3,0)-- (3,0);
        \draw [line width=1pt] (0,4.46)-- (0,-1);
        \filldraw  (-2,0) circle (3pt);
        \filldraw (0,3.46) circle (3pt);
        \draw (2.5,3.5) node {$\mathcal B_5$};
    \end{scope}
    \begin{scope}[xshift=80,yshift=-160]
        \draw [line width=1pt] (-2.5,-0.87)-- (0.5,4.33);
        \draw [line width=1pt] (2.5,-0.87)-- (-0.5,4.33);
        \draw [line width=1pt] (-3,0)-- (3,0);
        \draw [line width=1pt] (0,4.46)-- (0,-1);
        \filldraw  (-2,0) circle (3pt);
        \filldraw  (0,0) circle (3pt);
        \filldraw (0,3.46) circle (3pt);
        \draw (2.5,3.5) node {$\mathcal B_6$};
    \end{scope}
    \begin{scope}[xshift=160,yshift=-160]
        \draw [line width=1pt] (-2.5,-0.87)-- (0.5,4.33);
        \draw [line width=1pt] (2.5,-0.87)-- (-0.5,4.33);
        \draw [line width=1pt] (-3,0)-- (3,0);
        \draw [line width=1pt] (0,4.46)-- (0,-1);
        \filldraw (-2,0) circle (3pt);
        \filldraw (0,0) circle (3pt);
        \filldraw (2,0) circle (3pt);
        \filldraw (0,3.46) circle (3pt);
        \draw (2.5,3.5) node {$\mathcal B_7$};
    \end{scope}

    \begin{scope}[xshift=240,yshift=-160]
        \draw [line width=1pt] (-1.5, -0.6)-- (1.5, 4.6);
        \draw [line width=1pt] (-1.5,4.6)-- (1.5,-0.6);
        \draw [line width=1pt] (-3,2)-- (3,2);
        \draw [line width=1pt] (0,5)-- (0,-1);
        \filldraw (0,2) circle (3pt);
        \draw (2.5,3.5) node {$\mathcal C_1$};
    \end{scope}
\end{tikzpicture}
\caption{\small Possible configurations of lines on a completely reducible plane with four distinct lines and without a triangle. Singular points are marked with a bullet.}
\label{fig:no-triangle-conf-distinct}
\end{figure}
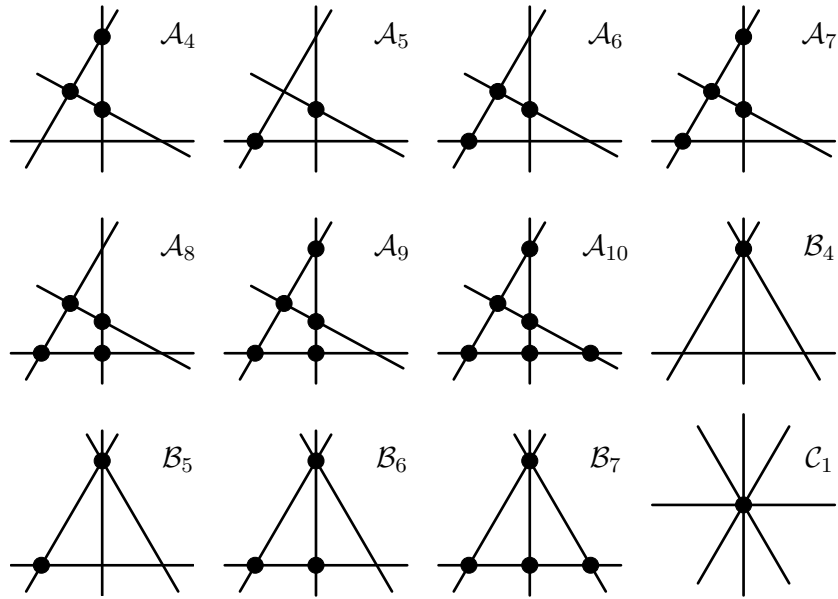

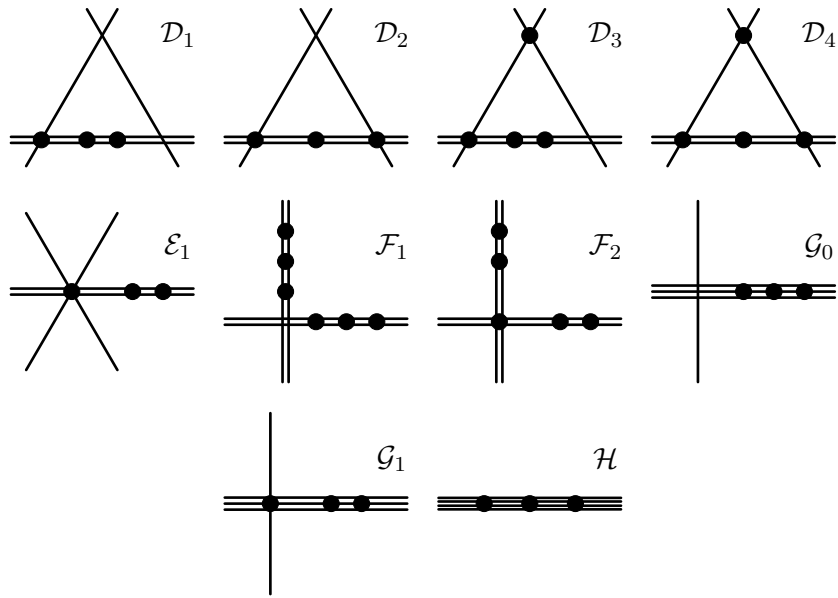
\begin{figure}
\vspace{0.5cm}
\centering
\begin{tikzpicture}[line cap=round,line join=round,>=triangle 45,x=0.4cm,y=0.4cm]
    \begin{scope}[xshift=0,yshift=0]
        \draw [line width=1pt] (-2.5,-0.87)-- (0.5,4.33);
        \draw [line width=1pt] (2.5,-0.87)-- (-0.5,4.33);
        \draw [line width=1pt] (-3,-0.1)-- (3,-0.1);
        \draw [line width=1pt] (-3,0.1)--(3,0.1);
        \draw (2.5,3.5) node {$\mathcal D_1$};
        \filldraw  (-2,0) circle (3pt);
        \filldraw  (-0.5,0) circle (3pt);
        \filldraw  (0.5,0) circle (3pt);
    \end{scope}
    \begin{scope}[xshift=80,yshift=0]
        \draw [line width=1pt] (-2.5,-0.87)-- (0.5,4.33);
        \draw [line width=1pt] (2.5,-0.87)-- (-0.5,4.33);
        \draw [line width=1pt] (-3,-0.1)-- (3,-0.1);
        \draw [line width=1pt] (-3,0.1)--(3,0.1);
        \draw (2.5,3.5) node {$\mathcal D_2$};
        \filldraw  (-2,0) circle (3pt);
        \filldraw  (0,0) circle (3pt);
        \filldraw  (2,0) circle (3pt);
    \end{scope}
    \begin{scope}[xshift=160,yshift=0]
        \draw [line width=1pt] (-2.5,-0.87)-- (0.5,4.33);
        \draw [line width=1pt] (2.5,-0.87)-- (-0.5,4.33);
        \draw [line width=1pt] (-3,-0.1)-- (3,-0.1);
        \draw [line width=1pt] (-3,0.1)--(3,0.1);
        \draw (2.5,3.5) node {$\mathcal D_3$};
        \filldraw  (-2,0) circle (3pt);
        \filldraw  (-0.5,0) circle (3pt);
        \filldraw  (0.5,0) circle (3pt);
        \filldraw (0,3.46) circle (3pt);
    \end{scope}
    \begin{scope}[xshift=240,yshift=0]
        \draw [line width=1pt] (-2.5,-0.87)-- (0.5,4.33);
        \draw [line width=1pt] (2.5,-0.87)-- (-0.5,4.33);
        \draw [line width=1pt] (-3,-0.1)-- (3,-0.1);
        \draw [line width=1pt] (-3,0.1)--(3,0.1);
        \draw (2.5,3.5) node {$\mathcal D_4$};
        \filldraw  (-2,0) circle (3pt);
        \filldraw  (0,0) circle (3pt);
        \filldraw  (2,0) circle (3pt);
        \filldraw (0,3.46) circle (3pt);
    \end{scope}
    \begin{scope}[xshift=0,yshift=-80]
        \draw [line width=1pt] (-2.5, -0.6)-- (0.5, 4.6);
        \draw [line width=1pt] (-2.5,4.6)-- (0.5,-0.6);
        \draw [line width=1pt] (-3,1.9)-- (3,1.9);
        \draw [line width=1pt] (-3,2.1)-- (3,2.1);
        \draw (2.5,3.5) node {$\mathcal E_1$};
        \filldraw  (-1,2) circle (3pt);
        \filldraw  (1,2) circle (3pt);
        \filldraw  (2,2) circle (3pt);
    \end{scope}
    \begin{scope}[xshift=80,yshift=-80]
        \draw [line width=1pt] (-0.9,5)-- (-0.9,-1);
        \draw [line width=1pt] (-1.1,5)-- (-1.1,-1);
        \draw [line width=1pt] (-3,0.9)-- (3,0.9);
        \draw [line width=1pt] (-3,1.1)-- (3,1.1);
        \draw (2.5,3.5) node {$\mathcal F_1$};
        \filldraw  (0,1) circle (3pt);
        \filldraw  (1,1) circle (3pt);
        \filldraw  (2,1) circle (3pt);
        \filldraw  (-1,2) circle (3pt);
        \filldraw  (-1,3) circle (3pt);
        \filldraw  (-1,4) circle (3pt);
    \end{scope}
    \begin{scope}[xshift=160,yshift=-80]
        \draw [line width=1pt] (-0.9,5)-- (-0.9,-1);
        \draw [line width=1pt] (-1.1,5)-- (-1.1,-1);
        \draw [line width=1pt] (-3,0.9)-- (3,0.9);
        \draw [line width=1pt] (-3,1.1)-- (3,1.1);
        \draw (2.5,3.5) node {$\mathcal F_2$};
        \filldraw  (-1,1) circle (3pt);
        \filldraw  (1,1) circle (3pt);
        \filldraw  (2,1) circle (3pt);
        \filldraw  (-1,3) circle (3pt);
        \filldraw  (-1,4) circle (3pt);
    \end{scope}
    \begin{scope}[xshift=240,yshift=-80]
        \draw [line width=1pt] (-1.5,5)-- (-1.5,-1);
        \draw [line width=1pt] (-3,1.8)-- (3,1.8);
        \draw [line width=1pt] (-3,2.0)-- (3,2.0);
        \draw [line width=1pt] (-3,2.2)-- (3,2.2);
        \draw (2.5,3.5) node {$\mathcal G_0$};
        \filldraw  (0,2) circle (3pt);
        \filldraw  (1,2) circle (3pt);
        \filldraw  (2,2) circle (3pt);
    \end{scope}
    \begin{scope}[xshift=80,yshift=-160]
        \draw [line width=1pt] (-1.5,5)-- (-1.5,-1);
        \draw [line width=1pt] (-3,1.8)-- (3,1.8);
        \draw [line width=1pt] (-3,2.0)-- (3,2.0);
        \draw [line width=1pt] (-3,2.2)-- (3,2.2);
        \draw (2.5,3.5) node {$\mathcal G_1$};
        \filldraw  (-1.5,2) circle (3pt);
        \filldraw  (0.5,2) circle (3pt);
        \filldraw  (1.5,2) circle (3pt);
    \end{scope}
    \begin{scope}[xshift=160,yshift=-160]
        \draw [line width=1pt] (-3,1.81)-- (3,1.81);
        \draw [line width=1pt] (-3,1.93)-- (3,1.93);
        \draw [line width=1pt] (-3,2.07)-- (3,2.07);
        \draw [line width=1pt] (-3,2.19)-- (3,2.19);
        \draw (2.5,3.5) node {$\mathcal H$};
        \filldraw  (-1.5,2) circle (3pt);
        \filldraw  (0,2) circle (3pt);
        \filldraw  (1.5,2) circle (3pt);
    \end{scope}
\end{tikzpicture}
\caption{\small Possible configurations of lines on a completely reducible plane with a multiple component and without a triangle. Singular points are marked with a bullet. Bullets on the same component, but not belonging also to another component can coincide.}
\label{fig:no-triangle-conf-double}
\end{figure}

\begin{lemma} \label{lem:no-triangle-conf}
If $X$ admits a completely reducible plane $\Pi$ without a triangle, then the intersection of $\Pi$ and $X$ has one of the configurations in Figure \ref{fig:no-triangle-conf-distinct}, if the lines on $\Pi$ are pairwise distinct, or Figure \ref{fig:no-triangle-conf-double}, if there is at least one multiple component.
\end{lemma}
\proof
We omit the proof, which employs the same combinatorial arguments as in Lemma~\ref{lem:triangle-conf}.
\endproof

The nomenclature for completely reducible planes introduced in Figures~\ref{fig:triangle-conf}, \ref{fig:no-triangle-conf-distinct} and \ref{fig:no-triangle-conf-double} will be used also in the next sections.

\vspace{\fill}

\section{Elliptic lines} \label{sec:char2-elliptic}

From now on, $\KK$ will denote an algebraically closed field of characteristic~2.

In this section we list some results on separable elliptic lines, especially the bounds on their valency, which are summarized in Table \ref{tab:char2-elliptic}. We postpone the study of inseparable lines (both elliptic and quasi-elliptic) to Lemma~\ref{lemma:ins=>cusp}.

\begin{table}
\caption{\small Known bounds for the valency of a separable elliptic line according to its kind, degree and singularity. Sharp bounds are marked with an asterisk *.}
\label{tab:char2-elliptic}
\vspace{0.2cm}

\centering
 \begin{tabular}{llll} \toprule
 kind & degree & singularity & valency \\
 \midrule
 \multirow{3}{*}{first kind}  & $3$ & 0         & $\leq 18$  \\
                              & $2$ & 1         & $\leq 13$  \\
                              & $1$ & 2 or 1    & $\leq  8$  \\
 \midrule
 \multirow{4}{*}{second kind} & $3$ & 0         & $\leq 20^*$ \\
                              & $2$ & 1         & $\leq 10$  \\
                              & $1$ & 2         & $\leq 9$  \\
                              & $1$ & 1         & $\leq 11$  \\
 \midrule
 --                           & $0$ & 3, 2 or 1 & $\leq 2^*$   \\
 \bottomrule
 \end{tabular}
\end{table}

Lines of the first kind, both elliptic and quasi-elliptic, have been treated in Proposition~\ref{prop:1stkind}.

The analysis of lines of the second kind in characteristic~2 is essentially the same as in characteristic~0. 
In fact, our arguments dealt with torsion sections of order~3 and fail only in characteristic~3.
The only discrepancy is that lines of degree~$3$ and $2$ have different ramification types, but this is compensated by wild ramification. We omit the proofs, but state the main facts that we will use later on.

\begin{lemma} \label{lemma:char2-deg3-ram}
    If $\ell$ has degree $3$, then it is separable and has ramification $2_4$, $2_2^2$, $2_2 3_2$ or $3_2^2$.
\end{lemma}
As in characteristic~0, the last ramification type in Lemma~\ref{lemma:char2-deg3-ram} deserves to be given a name.

\begin{definition} A line $\ell$ of degree $3$ is said to be \emph{special} if it is of the second kind and has ramification~$3_2^2$. 
\end{definition}

\begin{proposition} \label{prop:char2-degree3}
    If $\ell$ is a line of the second kind of degree $3$, then $v(\ell) \leq 20$; moreover, if $v(\ell) >16 $, then $\ell$ is special. If $v(\ell) = 19$, then the line has type $(p,q) = (6,1)$ and the 1-fiber is a ramified fiber of type $\I_n$, $n\geq 2$; if $v(\ell) = 20$, then the line has type $(p,q) = (6,2)$ and both 1-fibers are ramified fibers of type $\I_{n_1}$ and $\I_{n_2}$, $n_1,n_2\geq 2$.
\end{proposition}

We can also parametrize special lines as in \cite[Lemma 4.1]{veniani1} (see also \cite[Lemma 4.4]{char2}). Therefore, special lines induce a symplectic automorphism of degree~$3$ which permutes the lines in their 3-fibers (see \cite[Remark 4.2]{veniani1}).


\section{Quasi-elliptic lines} \label{sec:char2-qe}

In this section we study quasi-elliptic lines and find bounds on their valency. We summarize our results in Table \ref{tab:char2-qe}.

\begin{table}
\caption{\small Known bounds for the valency of a quasi-elliptic line according to its degree and singularity. Sharp bounds are marked with an asterisk.}
\vspace{0.2cm}
\centering
\label{tab:char2-qe}
    \begin{tabular}{cccl} 
      \toprule
      degree & singularity & & valency \\
      \midrule
      $3$ & 0         & & $\leq 16^*$ \\
      \multirow{2}{*}{$2$} & \multirow{2}{*}{$1$}         & cuspidal &  $\leq 19^*$ \\
       &          & not cuspidal &  $\leq 13$ \\
      $1$ & 2         & &  $\leq 8$ \\
      $1$ & 1         & &  $\leq 12$ \\
      $0$ & 3, 2 or 1 & &  $\leq 2$ \\
 \bottomrule
 \end{tabular}
\end{table}

We will first recall some general facts about quasi-elliptic fibrations in characteristic 2 (see \cite{bom-mum,rud-sha2}).

\begin{definition}
    The curve formed by the closure of the locus of singular
points on the irreducible fibers is called \emph{curve of cusps} or \emph{cuspidal curve}.
\end{definition}

The cuspidal curve of a quasi-elliptic fibration on a K3 surface $Z \rightarrow \PP^1$ (in particular, of a fibration induced by a quasi-elliptic line) is a smooth rational curve $K$ such that $K\cdot F = 2$. Only the fiber types 
\[
\II,\, \III,\, \I^*_{2n},\, \III^*,\, \II^*
\]
can appear in a quasi-elliptic fibration. We call such fibers \emph{quasi-elliptic fibers}. The restriction of the fibration $Z\rightarrow \PP^1$ to $K$ is an inseparable morphism of degree $2$. The cuspidal curve meets a reducible fiber in the following ways (multiple empty dots represent different possibilities):

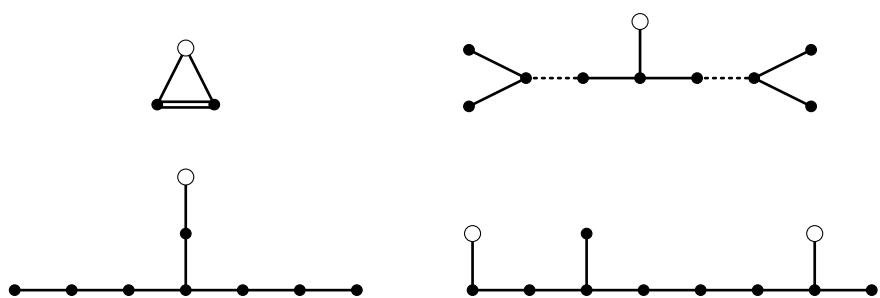
\begin{figure}[ht]
\centering
\begin{tikzpicture}[line cap=round,line join=round,>=triangle 45,x=0.75cm,y=0.75cm]
    \begin{scope}[xshift=0,yshift=-10]
        \filldraw (-0.5,0) circle (2pt);
        \filldraw (0.5,0) circle (2pt);
        \draw [line width=1pt] (-0.5,-0.05)-- (0.5,-0.05);
        \draw [line width=1pt] (-0.5, 0.05)-- (0.5, 0.05);

        \draw [line width=1pt] (0,1)--(-0.5,0);
        \draw [line width=1pt] (0,1)--(0.5,0);
        \filldraw [fill=white] (0,1) circle (3pt);
    \end{scope}

    \begin{scope}[xshift=170,yshift=0]
        \filldraw (-3,0.5) circle (2pt);
        \filldraw (-3,-0.5) circle (2pt);
        \filldraw (-2,0) circle (2pt);
        \filldraw (-1,0) circle (2pt);
        \filldraw (0,0) circle (2pt);
        \filldraw (1,0) circle (2pt);
        \filldraw (2,0) circle (2pt);
        \filldraw (3,0.5) circle (2pt);
        \filldraw (3,-0.5) circle (2pt);

        \draw [line width=1pt] (-3,0.5)--(-2,0)--(-3,-0.5);
        \draw [line width=1pt,dotted] (-2,0)--(-1,0);
        \draw [line width=1pt] (-1,0)-- (1,0);
        \draw [line width=1pt,dotted]  (1,0)-- (2,0);
        \draw [line width=1pt] (3,0.5)--(2,0)--(3,-0.5);

        \draw [line width=1pt] (0,0)--(0,1);
        \filldraw [fill=white] (0,1) circle (3pt);
    \end{scope}

    \begin{scope}[xshift=0,yshift=-80]
        \filldraw (-3,0) circle (2pt);
        \filldraw (-2,0) circle (2pt);
        \filldraw (-1,0) circle (2pt);
        \filldraw (0,0) circle (2pt);
        \filldraw (0,1) circle (2pt);
        \filldraw (1,0) circle (2pt);
        \filldraw (2,0) circle (2pt);
        \filldraw (3,0) circle (2pt);
        \draw [line width=1pt] (-3,0)-- (3,0);
        \draw [line width=1pt]  (0,0)-- (0,1);

        \draw [line width=1pt] (0,1)--(0,2);
        \filldraw [fill=white] (0,2) circle (3pt);
    \end{scope}

    \begin{scope}[xshift=150,yshift=-80]
        \filldraw (-2,0) circle (2pt);
        \filldraw (-1,0) circle (2pt);
        \filldraw (0,0) circle (2pt);
        \filldraw (0,1) circle (2pt);
        \filldraw (1,0) circle (2pt);
        \filldraw (2,0) circle (2pt);
        \filldraw (3,0) circle (2pt);
        \filldraw (4,0) circle (2pt);
        \filldraw (5,0) circle (2pt);
        \draw [line width=1pt] (-2,0)-- (5,0);
        \draw [line width=1pt]  (0,0)-- (0,1);

        \draw [line width=1pt] (-2,1)--(-2,0);
        \draw [line width=1pt]  (4,1)-- (4,0);
        \filldraw [fill=white] (-2,1) circle (3pt);
        \filldraw [fill=white]  (4,1) circle (3pt);
    \end{scope}
\end{tikzpicture}
\end{figure}

In particular, we observe that
\begin{itemize}
\item the way $K$ meets a reducible fiber is uniquely determined apart from type $\II^*$;
\item $K$ meets a component of a reducible fiber always transversally;
\item $K$ always meets only one component of multiplicity~$2$, with the exception of type $\III$, where it meets two components of multiplicity~$1$.
\end{itemize}

\begin{lemma} \label{lemma:sec-not-meet-cusp-curve}
A section of a quasi-elliptic fibration does not intersect the cuspidal curve.
\end{lemma}
\proof
We suppose the section and the cuspidal curve meet at a point on a fiber~$F$. Since the cuspidal curve intersects $F$ in a singular point or on a double component, the section would have intersection at least $2$ with $F$.
\endproof

\begin{lemma} \label{lemma:cusp-curve-not-exc-div}
The cuspidal curve of the fibration induced by a quasi-elliptic line $\ell$ cannot be an exceptional divisor.
\end{lemma}
\proof
Suppose the cuspidal curve of the fibration induced by $\ell$ coincides with an exceptional irreducible divisor~$E$; then $E$ must come from the resolution of a singular point $P$ on $\ell$. Since the general residual cubic is the image of a curve with a cusp on $E$, it should be singular in $P$, but this is ruled out by Lemma \ref{lemma:gen-res-cub-smooth}.
\endproof

Quasi-elliptic lines of degree $2$ play a special role. In fact, this is the only case where the strict transform $L$ of the line $\ell$ can be the cuspidal curve itself, since $L\cdot F = 2$. We give a name to these particular lines.

\begin{definition} \label{def:char2-cuspidal}
A line $\ell$ is said to be \emph{cuspidal} if it is quasi-elliptic of degree $2$ and the cuspidal curve on $Z$ coincides with the strict transform of $\ell$ itself.
\end{definition}

Assume that $\ell$ is a quasi-elliptic line which is \emph{not} cuspidal. By virtue of Lemma~\ref{lemma:cusp-curve-not-exc-div}, the cuspidal curve $K$ of the fibration induced by $\ell$ is a smooth rational curve in $Z$ of positive degree $k = K \cdot H>0$, where $H$ a hyperplane divisor in $Z$.

We introduce now a way to parametrize such cuspidal curves. Let us first choose a parameter $s$ for $K$ so that the restriction $\pi|_K:K\rightarrow \PP^1$ is given by 
\begin{equation} \label{eq:parameter-s=t^2}
    s \mapsto t = s^2,
\end{equation}
where $t$ parametrizes the planes containing $\ell$: $x_0 = tx_1$. Let $\bar K:= \rho(K)$ be the image of $K$ in $\PP^3$ through the resolution $\rho:Z\rightarrow X$. Then, we obtain a morphism from $K$ to $\bar K\subset \PP^3$ given by
\[
    \psi: s \mapsto \psi(s): = [\psi_0(s) : \psi_1(s) : \psi_2(s) : \psi_3(s) ],
\]
where $\psi_i(s)$ is a polynomial of degree $k$, $i=0,\ldots,3$.

Let now $\bar \pi: X\dashrightarrow \PP^1$ be the rational map $\bar\pi = \pi \circ \rho^{-1}$. Clearly, on the chart $x_1 \neq 0$, the map $\bar \pi$ can be written as 
\[
[x_0:x_1:x_2:x_3] \mapsto t = \frac{x_0}{x_1},
\]
since $\ell$ is not cuspidal and because of Lemma~\ref{lemma:cusp-curve-not-exc-div}, $\bar K$ has non-trivial intersection with the domain of definition of $\bar \pi$; hence, the restriction of $\pi = \bar \pi \circ \psi$ to (an open subset of) $K$ is given by
\[
    s \mapsto t = \frac{\psi_0(s)}{\psi_1(s)},
\]
Comparing with \eqref{eq:parameter-s=t^2}, we see that $\psi_0(s) = s^2\psi_1(s)$, i.e., we can parametrize $\bar K$ as
\begin{equation} \label{eq:K}
    s \mapsto \psi(s) = \left[ \, \sum_{i=0}^{k-2} a_i s^{i+2} : \sum_{i=0}^{k-2} a_i s^i: \sum_{i=0}^{k} b_i s^i : \sum_{i=0}^{k} c_i s^i \, \right],
\end{equation}
where $a_i,b_i,c_i \in \KK$, and $k$ is the degree of $K$. Note that the point $\psi(s)$ lies in the plane $x_0 = s^2 x_1$.

\subsection{Quasi-elliptic lines of degree 3}

Since a quasi-elliptic line of degree~$3$ is never cuspidal, in the following lemma we will suppose that the image of $K$ in $\PP^3$ is parametrized by~\eqref{eq:K}.

\begin{lemma} \label{lem:char2-d=3-degK}
If $\ell$ is a quasi-elliptic line of degree $3$, then the cuspidal curve of its fibration can have degree at most $4$.
\end{lemma}
\proof
The pullback of a plane containing $\ell$ is a divisor $H = F+L$ on the minimal resolution of $X$. The degree of $K$ is by definition $k := K \cdot H = K\cdot F + K\cdot L = 2 + K\cdot L$.

We claim that $K\cdot L\leq 2$. Indeed, $K$ cannot be the divisor $L$ itself, since $L\cdot F = 3$. If $K$ and $L$ meet at a point $P$, then the image of $P$ in $X$ (which we will call $P$ again) is a ramification point for $\ell$. Since there are at most two ramification points, $K\cdot L\leq 2$, unless $K$ is tangent to $L$ in at least one of them. Suppose therefore that $K$ is tangent to $L$ at $P$. 

\begin{claim}
The point $P$ cannot be of ramification $3$.
\end{claim}
\proof[Proof of the claim]
The fiber $F$ passing through $P$ must be singular in $P$, since $K$ is going through $P$. On the other hand, the corresponding residual cubic $C$ cannot be reducible, because $P$ is a point of ramification index $3$. 
In fact, if $C$ is the union of a line $m$ and an irreducible conic $Q$, then $Q$ is tangent to~$\ell$, and $m$ intersects $Q$ in two different points, giving rise to a fiber of type~$\I_n$: since $\ell$ is quasi-elliptic, this is not possible. 
If $C$ is the union of three lines, possibly not all distinct, then $K$ would have triple intersection with $F$, also impossible.
Therefore the cubic $F$ is irreducible, with a cusp in $P$ and tangent to $\ell$. Since $K$ is also tangent to $L$, we have $K\cdot F \geq 3$.
\endproof

We can thus assume that $P$ has ramification index $2$. Up to a change of coordinates, $P$ can be given by $[0:0:0:1]$. This means choosing $a_{0112} = 0$, $a_{0103} = 0$, while $a_{0121}$ and $a_{1003}$ must be non-zero, and can be normalized to~$1$.

\begin{claim}
The point $P$ is of ramification $2_4$.
\end{claim}
\proof[Proof of the claim]
Consider the parametrization~\eqref{eq:K}. Since $K$ goes through $P$, we can set $a_0 = 0$, $b_0 = 0$ and normalize $c_0 = 1$. Imposing that $\psi(s)$ is a singular point of the cubic in the plane $x_0 = s^2x_1$ for all $s$, one finds -- beside other relations -- that
\[
 a_{1012} = a_{0130},
\]
which is the condition for $P$ to be a point of ramification $2_4$.
\endproof

It follows that $\ell$ has only one ramification point, so that $K\cdot L \leq 2$ unless $K$ is tangent to $L$ of order at least $3$. Hence, we further set $a_2 = 0$, but this condition leads to $\ell$ containing a singular point, which contradicts the fact that $\ell$ has degree $3$.
\endproof

\begin{proposition} \label{prop:d=3-v}
If $\ell$ is a quasi-elliptic line of degree $3$, then $v(\ell) \leq 16$.
\end{proposition}
\proof
Let $iii$, $i_n^*$, $iii^*$, $ii^*$ be the numbers of reducible fibers of type $\III$, $\I_n^*$, $\III^*$, $\II^*$. The Euler--Poincaré characteristic formula (see \cite[Proposition 5.1.6]{cossec-dolgachev}) yields
\begin{equation} \label{eq:char2-euler-qe}
    iii + \sum (4+n)\,i_n^* + 7\,iii^* + 8\,ii^* = 20.
\end{equation}
A fiber of type $\III$ has local valency at most $1$ (because if it contains a line, the other component must be an irreducible conic), while other reducible fibers have valency at most $3$. 
Moreover, fibers of type $\III^*$ and $\II^*$ have valency $\leq 2$ because they do not contain three simple components. 
Hence, using equation~\eqref{eq:char2-euler-qe} we obtain a first rough estimate
\begin{equation} \label{eq:rough-estimate}
    v(\ell) \leq iii + 3\,\sum i_n^* + 2\,iii^* + 2\,ii^* \leq 20.
\end{equation}

We want to rule out all possible configurations that lead to $v(\ell) > 16$. We first note that a fiber of type $\III$ has local valency equal to either $1$ or $0$, according to whether the residual cubic is the union of a line and a conic, or an irreducible cubic with a cusp at a singular point of type $\bA_1$. We denote the number of the former $\III$-fibers with $iii'$ and of the latter with $iii''$.

First of all, if $v(\ell) >16$, then $iii^* = ii^* = 0$ and $i_n^* = 0$ for all $n > 2$, by formulas~\eqref{eq:char2-euler-qe} and \eqref{eq:rough-estimate}. If $i_2^* > 0$, then $i_2^* = 1$, $iii = iii'=14$, and the $\I_2^*$-fiber must have valency $3$; the corresponding cubic must then have the following shape:
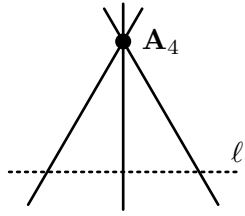
\begin{figure}[ht]
\centering
\begin{tikzpicture}[line cap=round,line join=round,>=triangle 45,x=0.5cm,y=0.5cm,
sing-pt/.style={circle,fill,inner sep=0pt,minimum size=6pt}]
    \begin{scope}[xshift=0,yshift=0]
        \draw [line width=1pt] (-2.5,-0.87)-- (0.5,4.33);
        \draw [line width=1pt] (2.5,-0.87)-- (-0.5,4.33);
        \draw [line width=1pt,dotted] (-3,0)-- (3,0);
            \node [above] at (3,0) {$\ell$};
        \draw [line width=1pt] (0,4.46)-- (0,-1);
        
        \node [sing-pt] at (0,3.46) [label=right:$\bA_4$] {};
    \end{scope}
\end{tikzpicture}
\end{figure}

In fact, the lines must meet at the same point (otherwise $\ell$ would have a fiber of type $\I_n$), and this point must be singular of type $\bA_4$ (the strict transforms of the three lines are three simple components of the fiber $\I_2^*$; the dual graph of the remaining components is an $\bA_4$-diagram).

If $i_2^*=0$ and $i_0^*>0$, then $i_0^*=1$ or $2$. Observe that there cannot be a fiber of type $\I_0^*$ not contributing to $v(\ell)$, otherwise $v(\ell) \leq 16$. Arguing as with type $\I_2^*$, it follows that the residual cubic of a fiber of type $\I_0^*$ is one of the types $\I^*_{0,n}$, $n=1,2,3$, as pictured:

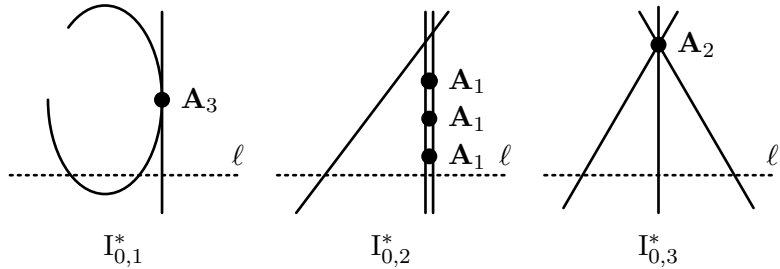
\begin{figure}[ht]
\centering
\begin{tikzpicture}[line cap=round,line join=round,>=triangle 45,x=0.5cm,y=0.5cm,
sing-pt/.style={circle,fill,inner sep=0pt,minimum size=6pt},]
    \begin{scope}[xshift=0,yshift=0]
        \draw [line width=1pt] (-2,2) arc (-180:130:1.5 and 2.5);
        \draw [line width=1pt,dotted] (-3,0)-- (3,0);
            \node [above] at (3,0) {$\ell$};
        \draw [line width=1pt] (1,4.33)-- (1,-1);
        \node [sing-pt] at (1,2) [label=right:$\bA_3$] {};
        \node at (0,-2) {$\I_{0,1}^*$};
    \end{scope}
    \begin{scope}[xshift=200,yshift=0]
        \draw [line width=1pt] (-2.5,-0.87)-- (0.5,4.33);
        \draw [line width=1pt] (2.5,-0.87)-- (-0.5,4.33);
        \draw [line width=1pt,dotted] (-3,0)-- (3,0);
            \node [above] at (3,0) {$\ell$};
        \draw [line width=1pt] (0,4.46)-- (0,-1);
        \node [sing-pt] at (0,3.46) [label=right:$\bA_2$] {};
        \node at (0,-2) {$\I_{0,3}^*$};
    \end{scope}
    \begin{scope}[xshift=100,yshift=0]
        \draw [line width=1pt] (1.5,4.33)-- (-2.5,-1);
        \draw [line width=1pt,dotted] (-3,0)-- (3,0);
            \node [above] at (3,0) {$\ell$};
        \draw [line width=1pt] (1.1,4.33)-- (1.1,-1);
        \draw [line width=1pt] (0.9,4.33)-- (0.9,-1);
        \filldraw  (1,2.5) circle (3pt);
        \node [sing-pt] at (1,2.5) [label=right:$\bA_1$] {};
        \node [sing-pt] at (1,1.5) [label=right:$\bA_1$] {};
        \node [sing-pt] at (1,0.5) [label=right:$\bA_1$] {};
        \node at (0,-2) {$\I_{0,2}^*$};
    \end{scope}
\end{tikzpicture}
\end{figure}

If $i_0^* = 2$, then $iii = 12$, one of the $\I_0^*$ must have valency~$3$ (i.e., be of type~$\I_{0,3}^*$) and the other must have valency~$3$ or $2$ (i.e., be of type $\I_{0,3}^*$ or $\I_{0,2}^*$). Once the numbers $i_2^*$ and $i_0^*$ are fixed, the sum $iii = iii' + iii''$ is uniquely determined by the number of the other fiber types and it is then a matter of listing all possibilities for $iii'$ and $iii''$ that lead to $v(\ell)>16$. There are 14 cases in total, as displayed in Table \ref{tab:d=3-v}.

\begin{table}[ht]
\caption{Fiber configurations for a quasi-elliptic line $\ell$ of degree $3$ with $v(\ell)>16$.}
\vspace{0.2cm}

\centering
\label{tab:d=3-v}
    \begin{tabular}{cccccccc} 
      \toprule
      case &$i_2^*$ & $i_{0,3}^*$ & $i_{0,2}^*$ & $i_{0,1}^*$ & $iii'$ & $iii''$ & $v(\ell)$ \\
      \midrule
      1&$1$     &             &             &             &   $14$ &         & $17$ \\
      2 &        &   $2$       &             &             &   $12$ &         & $18$ \\
      3 &        &   $2$       &             &             &   $11$ & $1$     & $17$ \\
      4 &        &   $1$       & $1$         &             &   $12$ &         & $17$ \\
      5 &        &   $1$       &             &             &   $16$ &         & $19$ \\
      6 &        &   $1$       &             &             &   $15$ & $1$     & $18$ \\
      7 &        &   $1$       &             &             &   $14$ & $2$     & $17$ \\
      8 &        &             & $1$         &             &   $16$ &         & $18$ \\
      9 &        &             & $1$         &             &   $15$ & $1$     & $17$ \\
      10&        &             &             &  $1$        &   $16$ &         & $17$ \\
      11&        &             &             &             &   $20$ &         & $20$ \\
      12&        &             &             &             &   $19$ & $1$     & $19$ \\
      13&        &             &             &             &   $18$ & $2$     & $18$ \\
      14 &        &             &             &             &   $17$ & $3$     & $17$ \\
 \bottomrule
 \end{tabular}
\end{table}

For each case, one can check that the lattice generated by $L$, the components of its fibers, the cuspidal curve $K$ (which has degree $k$, $2\leq k \leq 4$, according to Lemma \ref{lem:char2-d=3-degK}), and a general fiber has rank~$23$. 
\endproof

\begin{remark} \label{rmk:smooth=>ell-lines}
    An immediate corollary of the last proposition is the fact that on a smooth surface all lines are elliptic, which has already been proven by Rams and Schütt \cite[Proposition 2.1]{char2} using a different approach. In fact, a quasi-elliptic line on a smooth surface could only have 20 fibers of type $\III$, falling into case~11 of Table~\ref{tab:d=3-v} (which would be the only case to be ruled out).
\end{remark}

\begin{example}
The bound of Proposition \ref{prop:d=3-v} is sharp and is reached, for example, by the line $\ell:x_0 = x_1 = 0$ in the surface
\[
X: x_{0}^{3} x_{1} + x_{0} x_{1}^{3} + x_{0} x_{2}^{3} + x_{0}^{2} x_{1} x_{3} + x_{1} x_{2}^{2} x_{3} + x_{0}^{2} x_{3}^{2} + x_{1}^{2} x_{3}^{2} + x_{0} x_{3}^{3} = 0.
\]
The 16 lines meeting $\ell$ are given by $x_0 = x_3 = 0$ and by $x_1 = s^4 x_0,\, x_0 = a x_2 + b x_3$, where $a = s/(s^4 + s + 1)$, $b = 1/s^2$ and $s$ is a root of 
\[
s^{15}+s^{12}+s^9+s^8+s^7+s^6+s^4+s^3+s^2+s+1.
\]
\end{example}

\subsection{Quasi-elliptic lines of degree 2} 

\begin{proposition} \label{prop:d=2-v}
If $\ell$ is a quasi-elliptic line of degree $2$, then $v(\ell) \leq 19$.
\end{proposition}
\proof Consider the quasi-elliptic fibration induced by $\ell$; let $i$ be the number of fibers of type $\III$ and $j$ the number of reducible fibers not of type $\III$. The residual cubic of a fiber of type $\III$ cannot contain more than one line, while the other fibers contribute at most $2$ to the valency of $\ell$, since $\ell$ has degree $2$. It follows that
\[
v(\ell) \leq i + 2\,j.
\]
Computing the Euler--Poincaré characteristic yields $i+4\,j \leq 20$, so $v(\ell) \leq 20$.

We claim that $v(\ell)$ cannot be exactly $20$. Indeed, if $v(\ell) = 20$, then $\ell$ has exactly $20$ fibers whose components are a line and an irreducible conic.
Moreover, by Lemma \ref{lemma:exc-div-sec-or-fib-comp} the singular point on $\ell$ must be of type $\bA_1$, giving one divisor $E_0$ in the resolution which is a section of the fibration, or else the other divisors would form an extra fiber. 
The line $L$, the $20$ lines meeting $L$, the divisor $E_0$ and the general fiber $F$ generate a lattice of rank $23$, which is impossible.
\endproof
\begin{example}
The bound of Proposition \ref{prop:d=2-v} is sharp and is reached, for example, by the cuspidal line $\ell:x_0=x_1=0$ in the surface
\[
X : x_0^4 + x_1 x_2^3 + x_1^3 x_3 + x_0 x_2 x_3^2 = 0.
\]
The surface contains exactly one singular point $P = [0:0:0:1]$ of type $\bA_2$. The 19 lines meeting $\ell$ are given by $x_0 = t x_1,\, x_1 = t^5 x_2 + t^{15} x_3$, with $t$ any 19th root of unity. The surface contains exactly 20 lines. The line $\ell$ being quasi-elliptic, the surface $X$ has Picard number 22.
\end{example}

Such high valencies can indeed be reached only by cuspidal lines. To prove this fact we need to find a bound on the degree $k$ of the cuspidal curve $K$. Up to projective equivalence, we can assume that the singular point on the line is $P=[0:0:0:1]$ and that the cubic in $x_0 = s^2 x_1$ passes through $P$ twice for $s = 0$, and through $[0:0:1:0]$ for $s = \infty$. This means setting the following coefficients equal to zero:
\[
a_{1003},\, a_{0103}; \, a_{0112},\, a_{1030}.
\]
On the other hand, $a_{0130}$ and $a_{1012}$ must be non-zero in order to prevent $\ell$ from having degree 1; hence, we can normalize both of them to $1$.
In what follows, this will be the standard parametrization for lines of degree $2$.

The following lemma holds also for elliptic lines.

\begin{lemma} \label{lemma:ins=>cusp}
If $\ell$ is an inseparable line of valency $v(\ell)>12$, then $\ell$ is cuspidal.
\end{lemma}
\proof The line $\ell$ is inseparable if and only if $a_{1021} = a_{0121} = 0$. The smooth point of intersection of the residual cubic $E_t$ in $x_0 = s^2x_1$ with $\ell$ is given by
\[
P_t = [0:0:s:1].
\]
One can see explicitly that $P_s$ is a singular point of $E_s$ if and only if $s$ is a root of the following degree $6$ polynomial:
\begin{equation} \label{eq:char2-cusp-phi}
\begin{split}
        \varphi(s) & := a_{2020} s^{6} + a_{2011} s^{5} + a_{1120} s^{4} + a_{2002} s^{4} \\ & + a_{1111} s^{3} + a_{0220} s^{2} + a_{1102} s^{2} + a_{0211} s + a_{0202}.
\end{split}
\end{equation}
Furthermore, it can be checked by a local computation that if $E_s$ splits off a line, then $s$ is a root of $\varphi(s)$. Since there can be at most 2 lines through $P_s$, this implies that the valency of $\ell$ is not greater than $2\cdot 6 = 12$, unless the polynomial $\varphi$ vanishes identically, but $\varphi \equiv 0$ implies that all points $P_s$ are singular for $E_s$, i.e., the line $\ell$ is cuspidal.
\endproof
 

Let $\ell$ be a quasi-elliptic line of degree $2$, $K$ the cuspidal curve of the fibration induced by $\ell$, and $k = K\cdot H$ the degree of $K$. Considering the pullback of a plane containing $\ell$, we can write $H = F+L+E$, where $E = \sum n_i E_i$ has support on the exceptional divisors coming from the singular point on $\ell$. We will denote by $P$ the point of intersection of $E$ and $\ell$. The following lemma will help us determine the coefficients $n_i$.

\begin{lemma} \label{lemma:d=2-cubics}
If $\ell $ is a quasi-elliptic line of degree $d = 2$ and valency $v > 10$, then it is contained in a plane with a configuration~$\mathcal R_1$ or a configuration~$\mathcal R_2$.
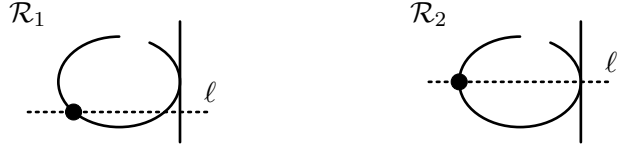
\begin{figure}[h!]
\centering
\begin{tikzpicture}[line cap=round,line join=round,>=triangle 45,x=0.4cm,y=0.4cm]
    \begin{scope}[xshift=0,yshift=0]
        \draw [line width=1pt] (0,2.5) arc (-270:60:2 and 1.5);
        \draw [line width=1pt, dotted] (-3,0)-- (3,0);
            \node [above] at (3,0) {$\ell$};
        \draw [line width=1pt] (2,-1)-- (2,3);
        \filldraw  (-1.5,0) circle (3pt);
        \node [above] at (-3,2.5) {$\mathcal R_1$};
    \end{scope}
    \begin{scope}[xshift=150,yshift=0]
        \draw [line width=1pt] (0,2.5) arc (-270:60:2 and 1.5);
        \draw [line width=1pt, dotted] (-3,1)-- (3,1);
            \node [above] at (3,1) {$\ell$};
        \draw [line width=1pt] (2,-1)-- (2,3);
        \filldraw (-2,1) circle (3pt);
        \node [above] at (-3,2.5) {$\mathcal R_2$};
    \end{scope}
\end{tikzpicture}
\end{figure}
\end{lemma}

\proof
Suppose there are no residual cubics as in the picture. Assume that $F$ is a fiber of $\ell$ of type $\III$ with $v(F) > 0$. 
Its residual cubic must split into a line and an irreducible conic; since its residual cubic cannot be of type $\mathcal R_1$ or $\mathcal R_2$, it must have configuration $\mathcal R_3$ (which can appear only once, since the intersection number of the general residual cubic with $\ell$ at the singular point is equal to~1).

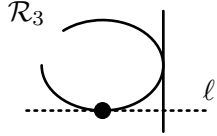
\begin{figure}[!ht]
\centering
\begin{tikzpicture}[line cap=round,line join=round,>=triangle 45,x=0.4cm,y=0.4cm]
    \draw [line width=1pt] (-2.5,1.5) arc (-180:130:2 and 1.5);
    \draw [line width=1pt, dotted] (-3,0)-- (3,0);
        \node [above] at (3,0) {$\ell$};
    \draw [line width=1pt] (1.5,3.3)-- (1.5,-0.7);
    \filldraw  (-0.5,0) circle (3pt);
    \node [above] at (-3,2.5) {$\mathcal R_3$};
\end{tikzpicture}
\end{figure}

On the other hand, reducible fibers $F$ not of type $\III$ with $v(F) > 0$ must have $e(F)\geq 6$ (because $\ell$ is quasi-elliptic) and $v(F) \leq 2$ (because $\ell$ has degree~$2$). Therefore, there can be at most $5$ of them ($4$, if configuration $\mathcal R_3$ appears) and $v(\ell)$ can be at most $10$.
\endproof

\begin{lemma} \label{lemma:d=2-deg-K}
Let $\ell$ be a separable quasi-elliptic line of degree $2$ contained in a plane with one of the residual cubics as in Lemma~\ref{lemma:d=2-cubics}. If $Q$ is the point of ramification on $\ell$, then only the following cases are possible:
\begin{enumerate}[(i)]
\item $k=2$, $K \cdot E = 0$, $Q \notin K$, $Q \neq P$;
\item $k=2$, $K \cdot E = 0$, $Q \notin K$, $Q = P$;
\item $k=3$, $K \cdot E = 0$, $Q \in K$, $Q \neq P$; 
\item $k=3$, $K \cdot E = 1$, $Q \notin K$, $Q \neq P$;
\item $k=4$, $K \cdot E = 1$, $Q\in K$, $Q \neq P$;
\item $k=4$, $K \cdot E = 1$, $Q\in K$, $Q = P$.
\end{enumerate}
\end{lemma}
\proof
Since $K$ is not a component of $H$, $K\cdot H \geq K \cdot F =2$. Considering a residual cubic as in Lemma \ref{lemma:d=2-cubics}, it is clear that the coefficients $n_i$ in $E$ must be equal to $1$, since the plane must correspond to a fiber of type $\I_n$ for the other line. Therefore we can write $E = E_0 + E_1 + \ldots + E_{n-1}$, where $E_0$ is a section and the other $E_i$'s are fiber components (necessarily of the same fiber). Therefore, $E\cdot K \leq 1$.

Moreover, by a local computation one can see that $K \cdot L \leq 1$. In fact, $K$ can intersect $\ell$ only in the point of ramification $Q$. The local computation is needed to rule out that $K$ might be tangent to $\ell$ in $Q$. It follows that $K \cdot L \leq 1$ and the only possible cases are those listed. In fact, if $P$ coincides with $Q$, then $K\cdot E =1$.
\endproof

\begin{proposition} \label{prop:v>13=>cusp}
If $\ell$ is a quasi-elliptic line of degree $2$ and valency $v(\ell) > 13$, then $\ell$ is cuspidal.
\end{proposition}
\proof
Thanks to Lemma \ref{lemma:ins=>cusp}, we can assume that $\ell$ is separable. We claim that a separable quasi-elliptic line of degree $2$ is always of the first kind, whence the bound on the valency follows.

In order to prove this, we parametrize such lines according to the possible values of the degree $k$ of the cuspidal curve $K$. Assume that the image of $K$ in $\PP^3$ is parametrized as in \eqref{eq:K}. By Lemma~\ref{lemma:d=2-cubics}, we can apply Lemma~\ref{lemma:d=2-deg-K}. According to the cases described there, we have the following conditions, after normalization:
\begin{enumerate}[(i)]
\item $k=2$, $a_{1021} = b_0 = c_0 = b_2 = c_2 = 0$, $a_0 = 1$;
\item $k=2$, $a_{0121} = b_0 = c_0 = b_2 = c_2 = 0$, $a_0 = 1$;
\item $k=3$, $a_{1021} = b_0 = c_0 = a_1 = c_3 = 0$, $a_0 = 1$, $b_3 \neq 0$; 
\item $k=3$, $a_{0121} = a_0 = b_0 = b_3 = c_3 = 0$, $a_1 = 1$, $c_0 \neq 0$;
\item $k=4$, $a_{1021} = a_0 = b_0 = a_2 = c_4 = 0$, $a_1 = 1$, $c_0\neq 0$, $b_4 \neq 0$;
\item $k=4$, $a_{0121} = a_0 = a_1 = b_0 = b_4 = c_4 = 0$, $a_2 = 1$, $c_0 \neq 0$.
\end{enumerate}
In fact, we can always choose $\psi(0)$ to be either $[0:0:0:1]$ or $[0:1:0:0]$, and $\psi(\infty)$ to be either $[0:0:1:0]$ or $[1:0:0:0]$.

We impose that $\psi(s)$ is indeed a singular point of the residual cubic in the plane $x_0 = s^2 x_1$ for every $s \in \PP^1$. It turns out that we can always express the following coefficients in terms of the others:
\begin{align*} a_{0211}, a_{1111}, a_{2011}; a_{0400}, a_{1300}, a_{2200}, a_{3100}, a_{4000}; \\ a_{0310}, a_{1210}, a_{2110}, a_{3010}; a_{0301}, a_{1201}, a_{2101}, a_{3001}. \end{align*}
(The first three turn out to be always equal to zero). In all cases one can verify that if in addition the conditions for being a line of the second kind are also satisfied, then all points on $K$ are singular, which contradicts the fact that $X$ only admits isolated singularities.
\endproof


\begin{example}
The following surface contains a separable quasi-elliptic line $\ell:x_0=x_1=0$ of degree $2$ and valency $12$:
\begin{multline*}
X : x_{0}^{3} x_{1} + x_{0}^{2} x_{1}^{2} + x_{0} x_{1}^{3} + x_{0}^{2} x_{1} x_{2} + x_{0} x_{1}^{2} x_{2} + x_{0}^{2} x_{2}^{2} \\
  + x_{0} x_{1} x_{2}^{2} + x_{1} x_{2}^{3} + x_{0} x_{1}^{2} x_{3} + x_{1} x_{2}^{2} x_{3} + x_{1}^{2} x_{3}^{2} + x_{0} x_{2} x_{3}^{2} = 0.
\end{multline*}
Two lines meeting $\ell$ are contained in the plane $x_1 = x_0$, while other 10 lines are contained in $x_1 = t x_0$, where $t$ is a root of
\[
t^{10} + t^{8} + t^{5} + t^{4} + 1.
\]
\end{example}

\subsection{Quasi-elliptic lines of degree 1}
In order to study quasi-elliptic lines of degree $d = 1$, we will need to find a bound on the degree $k = K\cdot H$ of the cuspidal curve $K$. Considering the pullback of a general plane containing $\ell$, one can write the hyperplane divisor as
\[
H = F + L + \sum n_i E_i
\]
where $F$ is a general fiber, $L$ is the strict transform of $\ell$ and the sum goes over the exceptional divisors coming from the singular points on $\ell$. Note that $L$ is a section of the quasi-elliptic fibration, so $K \cdot L = 0$. According to Lemma \ref{lemma:exc-div-sec-or-fib-comp}, for each singular point there is an exceptional divisor $E_0$ which is a section (hence $K\cdot E_0 = 0$), while the others are fiber components (hence $K \cdot E_i \leq 1$ and equality holds for at most two $E_i$'s per fiber). The following lemma will be useful to determine the coefficients $n_i$.

\begin{lemma} \label{lemma:d=1-cubics}
If $\ell$ is a line of degree $d = 1$ and valency $v \geq 5$ and singularity~$s$, then it is contained in a plane with one of the following residual cubics:
\end{lemma}

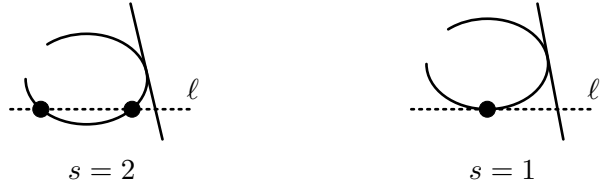
\begin{figure}[ht]
\centering
\begin{tikzpicture}[line cap=round,line join=round,>=triangle 45,x=0.4cm,y=0.4cm]
    \begin{scope}[xshift=0,yshift=0]
        \draw [line width=1pt] (-2.5,1) arc (-180:130:2 and 1.5);
        \draw [line width=1pt, dotted] (-3,0)-- (3,0);
            \node [above] at (3,0) {$\ell$};
        \draw [line width=1pt] (0.95,3.5)-- (2,-1);
        \filldraw  (1,0) circle (3pt);
        \filldraw  (-2,0) circle (3pt);
        \node at (0,-2) {$s = 2$};
    \end{scope}
    \begin{scope}[xshift=150,yshift=0]
        \draw [line width=1pt] (-2.5,1.5) arc (-180:130:2 and 1.5);
        \draw [line width=1pt, dotted] (-3,0)-- (3,0);
            \node [above] at (3,0) {$\ell$};
        \draw [line width=1pt] (1.15,3.5)-- (2,-1);
        \filldraw  (-0.5,0) circle (3pt);
        \node at (0,-2) {$s = 1$};
    \end{scope}
\end{tikzpicture}
\end{figure}

\proof
If $F$ is a fiber of type $\III$ with $v_\ell(F) = 1$ (the maximum possible since $d=1$), then its residual cubic contains a line meeting $\ell$ in a smooth point and an irreducible conic, i.e., it is one of the residual cubics as in the figure. Without fibers of type $\III$ and $v(F) =1$, $\ell$ cannot have valency greater than $5$: in fact, all other fibers either do not contribute to the valency of $\ell$ or have Euler--Poincaré characteristic $\geq 6$ and local valency~1.
\endproof

\begin{lemma} \label{lemma:d=1-deg-K}
If a quasi-elliptic line $\ell$ is contained in a plane with a residual cubic as in Lemma~\ref{lemma:d=1-cubics}, then its cuspidal curve has degree at least~$2$ and at most~$4$.
\end{lemma}
\proof Let us call $\ell'$ and $Q$ the line in the residual cubic on the plane~$\Pi$ given by Lemma~\ref{lemma:d=1-cubics}. 

Suppose first that $\ell$ has two singular points $P$ and $P'$. Let $H = F + L + \sum n_i E_i + \sum n_i'E_i'$ be the pullback of a general plane containing $\ell$, where $F$ is a general fiber of the fibration induced by $\ell$ and the $E_i$'s (resp. $E_i'$'s) come from the resolution of $P$ (resp. $P'$). We let $E_0$ and $E_0'$ be the sections of $\pi$ (Lemma~\ref{lemma:exc-div-sec-or-fib-comp}). 
Since the plane $\Pi$ corresponds to a fiber of type $\I_N$ for the line $\ell'$, the points $P$ and $P'$ must be of type $\bA_n$ and $\bA_{n'}$ respectively (with $N = n+n'+2$); moreover, pulling back $\Pi$ we see that the coefficients of the $E_i$'s and $E_i'$'s must be equal to $1$, so
\[ H = F+L +\sum_{i = 0}^n E_i + \sum_{i=0}^{n'} E_i', \] 
and we have the following diagram (curves of genus $0$ are marked with a circle, curves of genus $1$ with a square).

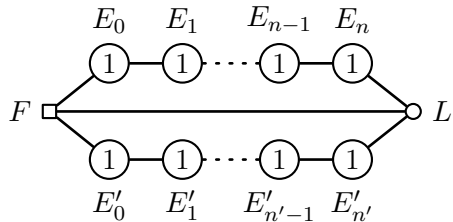
\begin{figure}[ht]
\centering
\begin{tikzpicture}[node distance = 12pt,
line cap=round,line join=round,>=triangle 45,x=0.75cm,y=0.75cm,
genus0/.style={circle,draw,thick,inner sep=2pt,minimum size=5pt},
genus1/.style={rectangle,draw,thick,inner sep=2pt,minimum size=5pt},]

    \coordinate (start);

    \node [genus1] (f)      [left=65pt of start,label=left:$F$]                 { };
    \node [genus0] (e_1)    [above left=18pt of start,label=above:$E_1$]        {1};
    \node [genus0] (e_n-1)  [above right=18pt of start,label=above:$E_{n-1}$]   {1};
    \node [genus0] (e_n-1') [below right=18pt of start,label=below:$E_{n'-1}'$] {1};
    \node [genus0] (e_1')   [below left=18pt of start,label=below:$E_1'$]       {1};
    \node [genus0] (e_0)    [left=of e_1,label=above:$E_0$]                     {1};
    \node [genus0] (e_n)    [right = of e_n-1,label=above:$E_n$]                {1};
    \node [genus0] (e_0')   [left = of e_1',label=below:$E_0'$]                 {1};
    \node [genus0] (e_n')   [right = of e_n-1',label=below:$E_{n'}'$]           {1};
    \node [genus0] (l)      [right=65pt of start,label=right:$L$]               { };
    
    \draw [line width=1pt] (f) -- (e_0) -- (e_1);
    \draw [line width=1pt,loosely dotted] (e_1) -- (e_n-1);
    \draw [line width=1pt] (e_n-1) -- (e_n) -- (l);
    \draw [line width=1pt] (f) -- (e_0') -- (e_1');
    \draw [line width=1pt,loosely dotted] (e_1') -- (e_n-1');
    \draw [line width=1pt] (e_n-1') -- (e_n') -- (l);
    \draw [line width=1pt] (f) -- (l);
\end{tikzpicture}
\end{figure}

The divisors $E_1,\ldots,E_{n-1}$ and $E_1',\ldots,E_{n-1}'$ are components of two distinct fibers for $\ell$, and the cuspidal curve $K$ can meet at most one of them in each fiber, so
\[
 K \cdot H = K \cdot F+ K\cdot L + K \cdot \sum E_i + K \cdot \sum E_i' \leq 2+0+1+1 = 4.
\]

Suppose now that $\ell$ has only one singular point $P$ (necessarily not of type~$\bA_1$, by explicit computation of the tangent cone). 

The pullback of $\Pi$ gives $H = F + L + \sum n_i E_i$.
The strict transforms $L$ and $\hat Q$ and the divisors $E_i$'s make up a fiber of the fibration induced by $\ell'$ which has at least two simple components ($L$ and $\hat Q$) and at least four components; hence, it must be of type $\I_{n+1}^*$, $\IV^*$ or $\III^*$. 
In the former case, we can have two different configurations, according to whether $P$ is of type $\bA_n$ or $\bD_n$, while in the latter two cases $P$ is of type $\bD_5$ resp. $\bE_6$. We have the following diagrams, where $E_0$ always denote the exceptional divisor which is a section for the fibration of~$\ell$ (Lemma~\ref{lemma:exc-div-sec-or-fib-comp}).

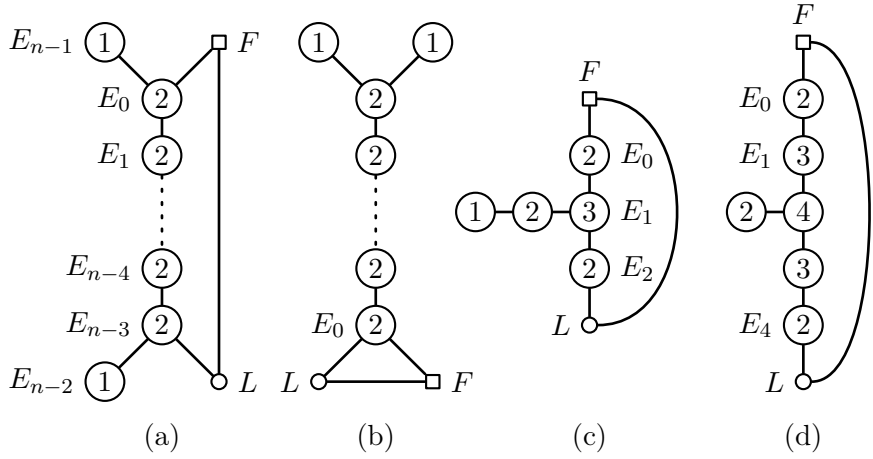
\begin{figure}[ht]
\centering
\begin{tikzpicture}[node distance = 12pt and 12pt,
line cap=round,line join=round,>=triangle 45,x=0.75cm,y=0.75cm,
genus0/.style={circle,draw,thick,inner sep=2pt,minimum size=5pt},
genus1/.style={rectangle,draw,thick,inner sep=2pt,minimum size=5pt},]

    \begin{scope}[xshift=0,yshift=0]
        \node at (0,-4) {(a)};
        
        \node [genus0] (e_n-1)  at (-1,3)  [label=left:$E_{n-1}$] {$1$};
        \node [genus1] (f)      at (1,3)   [label=right:$F$]                      {};
        \node [genus0] (e_0)    at (0,2)   [label=left:$E_{0}$]   {$2$};
        \node [genus0] (e_1)    at (0,1)   [label=left:$E_{1}$]   {$2$};
        \node [genus0] (e_n-4)  at (0,-1)  [label=left:$E_{n-4}$] {$2$};
        \node [genus0] (e_n-3)  at (0,-2)  [label=left:$E_{n-3}$] {$2$};
        \node [genus0] (e_n-2)  at (-1,-3) [label=left:$E_{n-2}$] {$1$};
        \node [genus0] (l)      at (1,-3)  [label=right:$L$]                      {};

        \draw [line width=1pt] (e_n-1)-- (e_0) -- (f) -- (l) -- (e_n-3) -- (e_n-2);
        \draw [line width=1pt] (e_0)-- (e_1);
        \draw [line width=1pt, loosely dotted] (e_1)-- (e_n-4);
        \draw [line width=1pt] (e_n-3)-- (e_n-4);
        
    \end{scope}
    
    \begin{scope}[xshift=80,yshift=0]
        \node [genus0] (e_n-1)  at (-1,3)    {$1$};
        \node [genus0] (e_n)    at (1,3)     {$1$};
        \node [genus0] (e_n-2)  at (0,2)     {$2$};
        \node [genus0] (e_n-3)  at (0,1)     {$2$};
        \node [genus0] (e_1)    at (0,-1)    {$2$};
        \node [genus0] (e_0)    at (0,-2)   [label=left:$E_0$] {$2$};
        
        \node [genus0] (l)    at (-1,-3)    [label=left:$L$] {};
        \node [genus1] (f)    at (1,-3)     [label=right:$F$] {};
    
        \draw [line width=1pt] (e_n-1) -- (e_n-2) -- (e_n);
        \draw [line width=1pt] (e_n-2) -- (e_n-3);
        \draw [line width=1pt, loosely dotted] (e_n-3) -- (e_1);
        \draw [line width=1pt] (e_1) -- (e_0) -- (l) -- (f) -- (e_0);
        
        \node at (0,-4) {(b)};
    \end{scope}
    
    \begin{scope}[xshift=160,yshift=0]
        \node [genus1] (f)      at (0,2)     [label=above:$F$] {};
        \node [genus0] (e_0)    at (0,1)     [label=right:$E_0$] {$2$};
        \node [genus0] (e_1)    at (0,0)     [label=right:$E_1$] {$3$};
        \node [genus0] (e_2)    at (0,-1)     [label=right:$E_2$] {$2$};
        \node [genus0] (e_3)    at (-1,0)     {$2$};
        \node [genus0] (e_4)    at (-2,0)     {$1$};
        \node [genus0] (l)      at (0,-2)     [label=left:$L$] {};

        \draw [line width=1pt] (f) -- (e_0) -- (e_1) -- (e_2) -- (l);
        \draw [line width=1pt] (e_4) -- (e_3) -- (e_1);
        \draw [line width=1pt] (f) .. controls (2,2) and (2,-2) .. (l);
        \node at (0,-4) {(c)};
    \end{scope}
    
    \begin{scope}[xshift=240,yshift=0]
        \node [genus1] (f)      at (0,3)     [label=above:$F$] {};
        \node [genus0] (e_0)    at (0,2)     [label=left:$E_0$] {$2$};
        \node [genus0] (e_1)    at (0,1)     [label=left:$E_1$] {$3$};
        \node [genus0] (e_2)    at (0,0)      {$4$};
        \node [genus0] (e_3)    at (0,-1)     {$3$};
        \node [genus0] (e_4)    at (0,-2)     [label=left:$E_4$] {$2$};
        \node [genus0] (e_5)    at (-1,0)      {$2$};
        \node [genus0] (l)      at (0,-3)     [label=left:$L$] {};
        
        \draw [line width=1pt] (f) -- (e_0) -- (e_1) -- (e_2) -- (e_3) -- (e_4) -- (l);
        \draw [line width=1pt] (e_5)-- (e_2);
        \draw [line width=1pt] (f) .. controls (1.5,3) and (1.5,-3) .. (l);
        \node at (0,-4) {(d)};
    \end{scope}
\end{tikzpicture}
\end{figure}

\begin{claim}
The cases (a) with $n\geq 5$, (c) and (d) are not possible.
\end{claim}
\proof[Proof of the claim] In these cases the exceptional divisors $E_1,\ldots,E_m$, with $m = n-2$, $4$ resp. $5$, are components of a fiber for $\ell$; we write this fiber as
\[
F = \sum_{i=1}^m r_i E_i + F',
\]
where $F'$ denotes the pullback of the residual cubic. Since $E_0$ and $L$ are sections, the $E_i$'s intersecting them must have multiplity $r_i =1$. But then we find a contradiction, since no quasi-elliptic fiber can have the following sub-configurations of divisors with multiplicities.

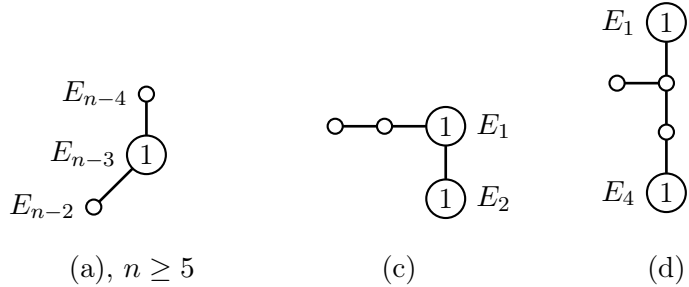
\begin{figure}[!ht]
\centering
\begin{tikzpicture}[node distance = 12pt and 12pt,
line cap=round,line join=round,>=triangle 45,x=0.75cm,y=0.75cm,
genus0/.style={circle,draw,thick,inner sep=2pt,minimum size=5pt},
genus1/.style={rectangle,draw,thick,inner sep=2pt,minimum size=5pt},]

    \begin{scope}[xshift=0,yshift=0]
        \coordinate [label=below:{(a), $n \geq 5$}] (start);
    
        \node [genus0] (e_n-2)  [above left=of start,label=left:$E_{n-2}$] { };
        \node [genus0] (e_n-3)  [above right=of e_n-2,  label=left:$E_{n-3}$] {1};
        \node [genus0] (e_n-4)  [above=of e_n-3,        label=left:$E_{n-4}$] { };

        \draw [line width=1pt] (e_n-4)-- (e_n-3) -- (e_n-2);
    \end{scope}
    
    \begin{scope}[xshift=100,yshift=0]
        \coordinate [label=below:(c)] (start);
        
        \node [genus0] (e_2)    [above right=of start,label=right:$E_2$]  {1};
        \node [genus0] (e_1)    [above=of e_2,  label=right:$E_1$]  {1};
        \node [genus0] (e_3)    [left=of e_1]                       { };
        \node [genus0] (e_4)    [left=of e_3]                       { };
        
        \draw [line width=1pt] (e_2) -- (e_1) -- (e_3) -- (e_4);
        
    \end{scope}
    
    \begin{scope}[xshift=200,yshift=0]
        \coordinate [label=below:(d)] (start);
        
        \node [genus0] (e_4) [above=of start,label=left:$E_4$]  {1};
        \node [genus0] (e_3) [above=of e_4]                     { };
        \node [genus0] (e_2) [above=of e_3]                     { };
        \node [genus0] (e_1) [above=of e_2,  label=left:$E_1$]  {1};
        \node [genus0] (e_5) [left =of e_2]                     { };

        \draw [line width=1pt] (e_4) -- (e_3) -- (e_2) -- (e_1);
        \draw [line width=1pt] (e_5)-- (e_2);
    \end{scope}
\end{tikzpicture}
\end{figure}
\endproof

In the few cases left, the multipliticities $n_i$ are not greater than $2$. 
\begin{enumerate}[(i)]
\item In case $P$ of type $\bD_n$, the $E_i$'s are part of the same fiber for $\ell$, so
\[ K \cdot \sum n_i E_i \leq 2.
\] 
\item In case $P$ of type $\bA_3$, $n_1 = n_2 = 1$, so $K \cdot (E_1 + E_2) \leq 2$. \item Finally, in case $P$ of type $\bA_4$, let
\[
 F = n_1 E_1 + n_2 E_2 + F'
\]
be the fiber containing $E_1$; since $E_1$ intersects the section $E_0$, $n_1 = 1$ and since $F$ has more than 3 components, $K$ does not intersect simple components, so $K\cdot E_1 = 0$. Thus,
\[
 K \cdot \sum n_i E_i = K\cdot (E_2 + E_3) \leq 2. \qedhere
\]
\end{enumerate}
\endproof

\begin{proposition}
If $\ell$ is a quasi-elliptic line of degree $1$ and singularity~$2$, then $v(\ell)\leq 8$.
\end{proposition}
\proof Up to coordinate change, we can suppose that the two singular points on the line $\ell$ are $P = [0:0:0:1]$ and $P' = [0:0:1:0]$. Moreover, we can assume that the residual cubic in $x_0 = t^2 x_1$ intersects $\ell$ twice in $P$ for $t = 0$, and twice in $P'$ for $t = \infty$. This means that the following coefficients can be set equal to zero:
\[ a_{1003},\, a_{0103};\, a_{1030},\, a_{0130}; \, a_{0112},\, a_{1021}; \]
whereas $a_{0121}$ and $a_{1012}$ must be non-zero and can be set equal to $1$.

Suppose that the image of the cuspidal curve $K$ in $\PP^3$ is parametrized by $\psi$ as in \eqref{eq:K}. By Lemmas \ref{lemma:d=1-cubics} and \ref{lemma:d=1-deg-K}, $K$ can have degree (i) $k = 2$, (ii) $k = 3$, or (iii) $k = 4$. 
It follows from the proof of the latter lemma that these cases happen exactly when the image of $K$ goes (i) neither through $P$ nor $P'$, (ii) through exactly one of them (say, $P$), or (iii) through both of them. 
When $\psi(s)$ is not equal to $P$ or $P'$ for $s=0$ or $s=\infty$, then up to a further coordinate change we can suppose that $\psi(0) = [0:1:0:0]$ or $\psi(\infty) = [1:0:0:0]$. 
Thus, we have the following conditions, after normalization:
\begin{enumerate}[(i)]
\item $k = 2$, $b_0 = c_0 = b_2 = c_2 = 0$, $a_0 = 1$;
\item $k = 3$, $a_0 = b_0 = b_3 = c_3 = 0$, $a_1 = 1$;
\item $k = 4$, $a_0 = b_0 = a_2 = c_4 = 0$, $a_1 = 1$.
\end{enumerate}

We then impose that the point $\psi(s)$ is in the zero locus of the derivatives of the residual cubics (parametrized by $x_0 = s^2 x_1$) for all $t \in \PP^1$, finding conditions on the coefficients $a_{i_0 i_1 i_2 i_3}$.

We observe that each residual cubic which contributes to the valency of $\ell$ must split off a line $m'$ passing through $\psi(s)$ and the point $[0:0:s^2:1]$ on the line $\ell$. 
The equations of $m'$ can be explicitly found (one of them is $x_0 = s^2x_1$ and the other is of the form $a x_1 + b x_2 + c x_3 = 0$) and imposing that $X$ contains $m'$ yields a polynomial in $s$, which is generically of degree $12-2\,k$ by an explicit computation.
This polynomial cannot be the zero polynomial, otherwise all points on the cuspidal curve would be singular, contrary to the fact that there are only isolated singularities on $X$. Therefore,
\[
v(\ell) \leq 12-2\,k \leq 8. \qedhere
\]
\endproof

\begin{proposition}
If $\ell$ is a quasi-elliptic line of degree $1$ and singularity~$1$, then $v(\ell)\leq 12$.
\end{proposition}
\proof
By Lemma \ref{lemma:d=1-cubics}, we can suppose that there exists a residual cubic splitting into a line $m$ and an irreducible conic tangent to $m$ and tangent to $\ell$ in its singular point $P$. Up to a change of coordinates, we can suppose that $P = [0:0:0:1]$, the line $m$ is given by $x_1 = x_3 = 0$, and the point of intersection of $m$ and $Q$ is $[1:0:0:0]$ (so that the plane $\Pi$ is given by $x_1 = 0$). Moreover, we impose that $\ell$ is of degree 1. All in all, this amounts to setting the following coefficients equal to zero in \eqref{eq:surfaceX}:
\[
a_{0103},\, a_{0112},\, a_{1003},\, a_{1012}; a_{4000},\, a_{3010},\, a_{2020},\, a_{1030};\, a_{3001},\, a_{2011}.
\]
On the other hand, to avoid contradictions such as $\ell$ having degree $0$ or the point $[1:0:0:0]$ being singular, the following coefficients must be non-zero:
\[ a_{1021},\, a_{0130};\, a_{3100}. \]
We then suppose that the degree of the cuspidal curve $K$ is $k$, so that the image of $K$ in $\mathbb P^3$ is parametrized by \eqref{eq:K}, and we impose that the coordinates of $\psi(s)$ satisfy the equations of the derivatives of the residual cubics relative to $\ell$ (parametrised by $x_0 = s^2 x_1$). We divide our analysis according to the value of $k$ which, by virtue of Lemma \ref{lemma:d=1-deg-K}, can be $2$, $3$ or $4$. When $k >2$ the curve $K$ must necessarily meet at least one exceptional divisor $E_i$ coming from $P$, so $\psi(t_0) = [0:0:0:1]$ for some $t_0$, and up to change of coordinates we can suppose that $t_0 = 0$. We can also normalize one of the $a_i$'s, once we know that it is non-zero. We divide the computations according to the following cases:
\begin{enumerate}[(i)]
\item $k = 2$, $a_0 = 1$;
\item $k = 3$, $a_0 = 0$, $a_1 = 1$;
\item $k = 4$, $a_0 = 0$, $a_1 = 1$;
\item $k = 4$, $a_0 = 0$, $a_1 = 0$, $a_2 = 1$.
\end{enumerate}

In each case, the choice of the following coefficients is unique (the first three turn out to be always equal to zero):
\begin{align*} a_{0211}, a_{1111}, a_{2011}; a_{0400}, a_{1300}, a_{2200}, a_{3100}, a_{4000}; \\ a_{0310}, a_{1210}, a_{2110}, a_{3010}; a_{0301}, a_{1201}, a_{2101}, a_{3001}. \end{align*}

We observe that if the residual cubic of a fiber contributing to the valency splits into a line $m'$ and an irreducible conic, then $m'$ passes through $\psi(s)$ and the point $[0:0:a_{1021}s^2:a_{0130}]$ on the line $\ell$. The equations of $m'$ can be explicitly found (one of them is $x_0 = s^2x_1$ and the other is of the form $ ax_1 + b x_2 + c x_3 = 0$) and imposing that $X$ contains $m'$ yields a polynomial in $s$, which is generically of degree $k+5$. This polynomial cannot be the zero polynomial, otherwise all points on the cuspidal curve would be singular, but there are only isolated singularities on $X$. Other residual cubics either do not contribute to the valency of $\ell$ or are singular in $P$, in which case they must be contained in the tangent cone of $P$ (since $P$ is not of type $\bA_1$, there can be at most two of them). Taking into account also the line $m$, we get that
\[
v(\ell) \leq (k+5) + 2 +1 \leq 12. \qedhere
\]
\endproof


\section{Proof of Theorem~\ref{thm:char2}} \label{sec:char2-proof}


\subsection{Triangle case}

\begin{proposition} \label{prop:triangle+sing-pt}
If $X$ contains a plane with a triangle and a singular point, then $\Phi(X) \leq 63$.
\end{proposition}
\proof
By Lemma \ref{lem:triangle-conf}, we can suppose that $X$ contains one of the configurations listed in Figure \ref{fig:triangle-conf}, except $\AAA_0$, $\BBB_0$ or $\CCC_0$ because they do not contain a singular point. 
In the configurations $\mathcal A_i$ and $\mathcal B_i$, $i=1,2,3$, all lines must be elliptic (because the corresponding fiber is of type $\I_n$ or $\IV$) and those of singularity $0$ cannot have valency higher than $18$ (since there cannot be an automorphism of degree $3$ exchanging the other three lines). Hence, for each line of the configuration we consider the following bounds on its valency $v$ according to its singularity $s$: if $s = 0$, $v \leq 18$; if $s = 1$, $v \leq 13$; if $s = 2$, $v\leq 11$; if $s = 3$, $v\leq 2$. Moreover, we use the fact that there can be at most $8$ lines going through each singular point (Lemma~\ref{lemma:linesthroughsingularpoint}). 
For configurations $\mathcal D_0$ and $\mathcal E_0$ we use the bound $v\leq 19$ for the simple lines, which might be quasi-elliptic.
We obtain the following bounds on the total number of lines $N$:
\begin{itemize}
\item configurations $\mathcal A_1$ and $\mathcal B_1$: $N\leq 62$;
\item configurations $\mathcal A_2$ and $\mathcal B_2$: $N\leq 63$;
\item configurations $\mathcal A_3$ and $\mathcal B_3$: $N\leq 57$;
\item configurations $\mathcal D_0$ and $\mathcal E_0$: $N\leq 58$. \qedhere
\end{itemize}
\endproof

We can now prove Theorem~\ref{thm:char2} in the triangle case.

\begin{proposition} \label{prop:char2-triangle}
If $X$ has a triangle, then $\Phi(X) \leq 68$.
\end{proposition}
\proof
By Lemma~\ref{lem:triangle-conf} and Proposition~\ref{prop:triangle+sing-pt} we can assume that $X$ contains a plane $\Pi$ with four distinct lines and no singular point (configurations $\AAA_0$, $\BBB_0$ or $\CCC_0$). If all lines have valency less or equal than $19$, then $X$ has at most $4\cdot(19-3)+4 = 68$ lines.

Suppose that $X$ has a line $\ell_0$ of valency $20$. Then $\ell_0$ is a special line of type $(6,2)$, by Tables \ref{tab:char2-elliptic} and \ref{tab:char2-qe} and Proposition \ref{prop:char2-degree3}, and the 1-fibers are the only ramified fibers. Moreover, the six 3-fibers do not contain any singular point, since if one of them did, then on account of the automorphism $\sigma$ induced by $\ell_0$ it would have Euler--Poincaré characteristic at least 6 and $e(Z)$ would exceed 24.

Let $\ell_i$, $i=1,\,2,\,3$, be the three other lines on $\Pi$ (which must be one of the unramified 3-fibers), let $m_1$ and $m_2$ be the lines in the 1-fibers, and $P_i$ the points on $m_i$ which sit on the residual conic but not on $\ell_0$ ($i=1,\,2$). 

The lines $\ell_i$ have the same valency, which must be greater than $18$, or else $X$ has at most $(20-3)+3\cdot(18-3)+4 = 66$ lines; moreover, they cannot be quasi-elliptic because they have an $\I_3$-fiber. Therefore, they induce an automorphism of degree $3$, whence they must have the same fibration of $\ell_0$; in particular, their 3-fibers do not contain singular points and they have valency~$20$.

\begin{claim}
The points $P_i$ are singular and the lines $m_i$ are cuspidal.
\end{claim}

\begin{figure}[ht]
\centering
\begin{tikzpicture}[line cap=round,line join=round,>=triangle 45,x=0.5cm,y=0.5cm,]
        \draw [line width=1pt] (-8,0)-- (8,0);
            \node [above] at (8,0) {$\ell_0$};
        \draw [line width=1pt] (-2.5,-0.87)-- (0.5,4.33);
        \draw [line width=1pt] (-2.13,2.23)-- (2.87,-0.5);
        \draw [line width=1pt] (0,4.46)-- (0,-1);
            \node at (-2.3,-1.5) {$\ell_1$};
            \node at ( 0,-1.5) {$\ell_2$};
            \node at ( 2.87,-1.5) {$\ell_3$};
            
            \node at (-2,3.75) {$6\times$};

        \draw [line width=1pt] (-6,-1)--(-6,4.33);
            \node at (-6,-1.5) {$m_1$};
        \draw [line width=1pt] (-7,1.75) arc (-180:130:1 and 1.75);
            \node [above left] at (-6,3.5) {$P_1$};
        \filldraw (-6,3.5) circle (3pt);
        
        \draw [line width=1pt] (6,-1)--(6,4.33);
            \node at ( 6,-1.5) {$m_2$};
        \draw [line width=1pt] (5,1.75) arc (-180:130:1 and 1.75);
            \node [above right] at (6,3.5) {$P_2$};
        \filldraw (6,3.5) circle (3pt);
\end{tikzpicture}
\end{figure}

\proof[Proof of the claim]
The lines $\ell_i$ also have fibration $(6,2)$. Let $n$ be a line in a 3-fiber of one of the $\ell_i$ which meets $m_1$: there are $15$ of them. By the same argument as for the $\ell_i$'s, $n$ must also be special lines of valency $20$. Note that, regardless of $P_1$ being singular or not, they must meet $m_1$ in a point different from $P_1$. 

It follows that $m_1$ must belong to a 1-fiber of $n$, because it cannot have valency $20$: therefore, $m_1$ has $16$ fibers of type $\III$ with ramification of order $2$. By the Riemann--Hurwitz formula, $m_1$ must be inseparable, so $P_1$ must be singular. Recalling Lemma~\ref{lemma:ins=>cusp}, we conclude that $m_1$ is cuspidal. The same reasoning applies to $P_2$ and $m_2$.
\endproof

The points $P_i$ must be of type $\bA_1$, otherwise the Euler--Poincaré characteristic of $X$ would be greater than $24$. Let now $C$ be the residual cubic contained in the plane on which both $m_1$ and $P_2$ lie. We claim that $C$ can be neither irreducible nor reducible, thus finding a contradiction.

In fact, if $C$ is irreducible, then $C$ must have a cusp in $P_2$, but this is impossible, since the cuspidal curve of the fibration induced by $m_1$ is the strict transform of $m_1$. On the other hand, if $C$ is reducible, then it has $2$ or $3$ components and, since $P_2$ is of type $\bA_1$, the corresponding fiber has $3$ or $4$ components, but there do not exist quasi-elliptic fibers with $3$ or $4$ components. \endproof


\subsection{Square case}

We employ here the technique of the line graph $\Gamma = \Gamma(X)$, which we recalled in §\ref{subsec:triangle-free}. Here, as $\Char \KK = 2$, a parabolic subgraph $D\subset \Gamma$ might induce an elliptic or a quasi-elliptic fibration. 

\begin{lemma} \label{lemma:D-ell}
If $D$ induces an elliptic fibration, then
\[
    \Phi(X) \leq v(D) + 24.
\]
\end{lemma}
\proof
The same proof as in \cite[Proposition 5.5]{veniani1} applies.
\endproof

\begin{lemma} \label{lem:square+v13=>68}
Let $X$ be a triangle-free K3 quartic surface $X$ with a square. If all lines have valency at most $13$, then $X$ contains at most $68$ lines.
\end{lemma}
\proof
The square $D$ induces an elliptic fibration because quasi-elliptic fibration cannot have fibers of type $\I_4$. Hence,
\[
    \Phi(X) \leq v(D) + 24 \leq 4\cdot(13-2) + 24 = 68.\qedhere
\]
\endproof

It is therefore important to classify all lines of valency greater than $13$ on triangle-free surfaces. By Lemma \ref{lem:elliptic-triangle-free}, all such lines must be quasi-elliptic.

\begin{proposition} \label{prop:comp-red-Pi=>leq68}
    If $X$ is a triangle free surface and admits a completely reducible plane, then $\Phi(X) \leq 68$.
\end{proposition}
\proof
The proof is a case-by-case analysis on the configurations of Figures \ref{fig:no-triangle-conf-distinct} and \ref{fig:no-triangle-conf-double} given by Lemma \ref{lem:no-triangle-conf}. We use the bound of Lemma \ref{lem:elliptic-triangle-free} for elliptic lines and those of Table \ref{tab:char2-qe}. Moreover, we use the fact that there are at most $8$ lines through a singular point (Lemma~\ref{lemma:linesthroughsingularpoint}).

We have to refine our argument only for configurations $\mathcal A_8$ and $\mathcal C_1$. 
\begin{itemize}
\item In configuration $\mathcal A_8$, all lines are elliptic, have degree $1$ and singularity $2$. Since the plane corresponds to an $\I_n$-fiber, with $n\geq 5$, they can be met by at most $9$ other lines in other planes. It follows that $X$ contains at most $4\cdot 9 + 4\cdot (8-2) + 4 = 64$ lines.
\item In configuration $\mathcal C_1$, if one of the lines is quasi-elliptic, then the plane corresponds to a fiber of type $\I_0^*$; by an Euler--Poincaré characteristic argument, the valency of the lines is not greater than $16$, so $X$ contains at most $4\cdot 16 + 4 = 68$ lines.
\end{itemize}

Note that $\mathcal C_1$ is the only configurations where $68$ can be reached.
\endproof

Now that we have ruled out completely reducible planes, it will be easier to classify lines with valency greater than 13. Since $X$ is triangle-free, such lines must be quasi-elliptic of degree $3$ or $2$.

\begin{lemma} \label{lemma:d=3-v>13}
If $\ell\subset X$ is a quasi-elliptic line of degree $3$ and valency $v>13$ on a surface without completely reducible planes, then it has one of the fibrations listed in Table~\ref{tab:d=3-13<v<=16-true}.
\end{lemma}
\proof
Since $X$ does not have completely reducible planes, a residual cubic of $\ell$ can be either irreducible or split into a line and a conic. In the former case, it must have a cusp, which might be a singular point of the surface; in the latter case, the line and the conic must be tangent, and their intersection point might be a singular point of the surface. We can have the following possibilities for reducible fibers (the extra subscript number denotes the local valency, while $a$ and $b$ distinguish the two possibilities for $\I^*_{2n,1}$):
\begin{description}[align=left,labelwidth=0.06\linewidth]
\item[$\III_0$]: cusp with a point of type $\bA_1$;
\item[$\III_1$]: line and conic with a smooth intersection point;
\item[$\I^*_{2n,0}$]: cusp with a point of type $\bD_{2n+4}$;
\item[$\I^{*a}_{2n,1}$]: line and conic with a point of type $\bA_{2n+3}$, $n \geq 0$;
\item[$\I^{*b}_{2n,1}$]: line and conic with a point of type $\bD_{2n+3}$, $n \geq 1$;
\item[$\III_0^*$]: cusp with a point of type $\bE_7$;
\item[$\III_1^*$]: line and conic with a point of type $\bE_6$;
\item[$\II^*_0$]: cusp with a point of type $\bE_8$.
\end{description}

We then make a list of the fibrations that lead to $v(\ell)>13$, imitating the arguments of Proposition \ref{prop:d=3-v}. The results are shown in Table \ref{tab:d=3-13<v<=16}. There are 14 cases; in two of them (cases 4 and 9) one has to distinguish the type of the fibers $\I^{*a}_{2n,1}$ and $\I^{*b}_{2n,1}$.

\begin{table}[t]
\caption{Candidates to the fiber configuration of a quasi-elliptic line $\ell$ of degree $3$ on a surface without completely reducible planes with $13 < v(\ell) \leq 16$.}
\label{tab:d=3-13<v<=16}
\vspace{0.2cm}
\centering
\begin{tabular}{c|ccccccc|c}
    \toprule
    case &$iii^*_1$ & $i_{2,1}^{*a,b}$ & $i_{2,0}^*$ & $i_{0,1}^*$ & $i_{0,0}^*$ & $iii_1$ & $iii_0$ & $v(\ell)$ \\
    \midrule
    1           &     &   &   & 1 &   & 15 & 1 &    16 \\
    2           &     &   &   &   & 1 & 16 & 1 &    16 \\
    3           &     &   &   &   &   & 16 & 4 &    16 \\
    $4^{a,b}$   &   & 1 &   &   &   & 14 &   &    15 \\
    5           &     &   &   & 1 &   & 14 & 2 &    15 \\
    6           &     &   &   &   & 1 & 15 & 1 &    15 \\
    7           &     &   &   &   &   & 15 & 5 &    15 \\
    8           & 1   &   &   &   &   & 13 &   &    14 \\
    $9^{a,b}$   &   & 1 &   &   &   & 13 & 1 &    14 \\
    10          &     &   & 1 &   &   & 14 &   &    14 \\
    11          &     &   &   & 2 &   & 12 &   &    14 \\
    12          &     &   &   & 1 &   & 13 & 3 &    14 \\
    13          &     &   &   &   & 1 & 14 & 2 &    14 \\
    14          &     &   &   &   &   & 14 & 6 &    14 \\
 \bottomrule
 \end{tabular}
\end{table}

We then explicitly compute the rank of the intersection matrix of these 14 cases, taking into account the general fiber $F$, the line $L$, the fiber components and the cuspidal curve $K$. The only unknown intersection number is $K\cdot L$, but by Lemma~\ref{lem:char2-d=3-degK} this can only be $0$, $1$ or $2$. If the rank is always greater than~$22$, the fibration is discarded. The cases that pass this test are those listed in Table~\ref{tab:d=3-13<v<=16-true} (all of them with $K\cdot L = 0$).
\endproof

\begin{table}[t]
\caption{Fiber configurations from Table \ref{tab:d=3-13<v<=16} for a quasi-elliptic line $\ell$ of degree $3$ with $13 < v(\ell) \leq 16$ which generate a lattice of rank $\leq 22$.}
\label{tab:d=3-13<v<=16-true}
\vspace{0.2cm}
\centering
\begin{tabular}{c|ccccccc|cc}
    \toprule
    case    &$iii^*_1$ & $i_{2,1}^*$ & $i_{0,1}^*$ & $i_{0,0}^*$ & $iii_1$ & $iii_0$ & $v(\ell)$ & $\Sing(X)$\\
    \midrule
    2       &    &   &   & 1 & 16 & 1 &    16   & $\bD_4, \bA_1$ \\
    3       &    &   &   &   & 16 & 4 &    16   & $4 \bA_1$\\
    $4^{b}$ &    & 1 &   &   & 14 &   &    15   & $\bD_5$ \\
    5       &    &   & 1 &   & 14 & 2 &    15   & $\bA_3, 2 \bA_1$\\
    8       & 1  &   &   &   & 13 &   &    14   & $\bE_6$ \\
    $9^{a}$ &    & 1 &   &   & 13 & 1 &    14   & $\bA_5, \bA_1$ \\
    11      &    &   & 2 &   & 12 &   &    14   & $2 \bA_3$ \\
    \bottomrule
\end{tabular}
\end{table}

\begin{lemma} \label{lemma:d=2-v>13}
If $\ell$ is a quasi-elliptic line of degree~$2$ and valency $v>13$ on a surface without completely reducible planes, then it has exactly $19$ fibers of type~$\III$ with valency~$1$ and the surface contains only one singular point of type~$\bA_2$.
\end{lemma}
\proof
By Lemma~\ref{lemma:ins=>cusp}, $\ell$ is cuspidal and we can assume that it admits the second configuration of Lemma~\ref{lemma:d=2-cubics}. Therefore, the singular point $P$ on $\ell$ is of type $\bA_n$. Let $F_0$ be the only fiber whose residual cubic intersects $\ell$ only in $P$, and let $F$ be a reducible fiber different from $F_0$. The residual cubic of $F$ cannot be irreducible (because its cusp is on $\ell$, which does not have other singular points, and thus $F$ would be of type $\II$), and it cannot split into three lines by hypothesis. It follows that $F$ is the union of a line and an irreducible conic, necessarily meeting tangentially at a point of $\ell$; hence, $F$ is of type $\III$ and has valency~$1$. 

On the other hand, the fiber $F_0$ is a reducible fiber: in fact, if it were not, then $\ell$ would have $20$ fibers of type $\III$ and valency $20$, but this is excluded by Proposition~\ref{prop:d=2-v}. 
The residual cubic $C_0$ of $F_0$ is irreducible (necessarily, with a cusp in $P$): indeed, it can never be the union of an irreducible conic and a line, because it would result in a fiber of type $\I_n$, and by hypothesis it cannot split into three lines. 
It follows that the strict transform of $C_0$ is a simple component of $F_0$, and the remaining components are supplied by the exceptional divisors coming from $P$. 
Since $P$ is of type~$\bA_n$, this is only possible if $F$ is of type~$\III$. 
By Lemma~\ref{lemma:exc-div-sec-or-fib-comp}, $P$ is of type $\bA_2$ (one exceptional irreducible divisor is a fiber component, the other one a section) and an Euler--Poincaré characteristic argument yields that there are $19$ other $\III$-fibers.
\endproof

\begin{proposition} \label{prop:v>13=>leq64}
Let $X$ be a K3 quartic surface without completely reducible planes. If $X$ contains a line $\ell$ of valency $v > 13$, then $\Phi(X) \leq 64$.
\end{proposition}
\proof
The proof is to be done case by case according to the possible fibrations of $\ell$ given by Lemma \ref{lemma:d=3-v>13} and Lemma \ref{lemma:d=2-v>13}. The fibrations univocally determine the singular locus of $X$ (for $d = 3$, we refer to Table \ref{tab:d=3-13<v<=16-true}) and are mutually exclusive, in the sense that the same surface cannot have two lines of valency greater than $13$ with different fibrations. We do one case as example; one can argue analogously for the other cases.

Suppose that $X$ contains a line falling in case $2$ of Table \ref{tab:d=3-13<v<=16-true}. The surface $X$ has then two singular points of type $\bD_4$ and $\bA_1$, $\ell$ is of degree~$3$ and $v(\ell) = 16$. Consider the set $S(\ell)$ of lines that do not meet $\ell$ (and are therefore sections of its fibration). The number of lines on $X$ is not greater than
\[
\# S(\ell) + v(\ell) + 1 = \# S(\ell) + 17.
\]
If all lines in $S(\ell)$ pass through the singular points, then there can be at most $16$ of them. Suppose then that $s\in S(\ell)$ does not go through the singular points. Then, by inspection of the intersection matrix, $s$ must meet exactly $8$ lines contained in the $\III_1$-fibers of $\ell$. We choose $2$ of these $8$ lines, say $m_1$ and $m_2$, such that the corresponding fibers are not ramified for $\ell$ (this is possible since ramification occurs in at most two fibers).

Now, $\ell$, $m_1$, $m_2$ and $s$ form a square $D$. The lines $m_1$ and $m_2$ must be elliptic, since they have a fiber of type $\I_2$. Since the valency of $s$ cannot be greater than 16, the number of lines on $X$ is not greater than
\begin{align*}
v(D) + 24 & = (v(m_1)-2)+(v(m_2)-2) + 8 + \\ & + (v(s)-2-8) + (v(\ell)-2-8) + 24 \\ & \leq 10+10+8+6+6+24 = 64. \qedhere
\end{align*}
\endproof

\begin{corollary} \label{cor:triangle-free+square}
    If $X$ is a triangle free surface containing a square, then $\Phi(X)\leq 68$.
\end{corollary}
\proof
If $X$ admits a completely reducible plane, then we can use Proposition~\ref{prop:comp-red-Pi=>leq68}. Otherwise, we conclude by Proposition \ref{prop:v>13=>leq64} and Lemma \ref{lem:square+v13=>68}.
\endproof

\subsection{Square-free case} 
We prove an analog of Lemma \ref{lemma:D-ell} for quasi-elliptic fibration.

\begin{lemma} \label{lemma:D-qe}
If $D\subset \Gamma$ induces a quasi-elliptic fibration, then 
\[
    \Phi(X) \leq v(D) + 25.
\]
\end{lemma}
\proof
We observe that a fiber of type $\III$ cannot contain two lines, since two lines never intersect tangentially. Therefore, applying formula \eqref{eq:char2-euler-qe} twice, one gets
\[
\begin{split}
\Phi(X) & \leq v(D) + iii + \sum (5+n)\,i_n^* + 8\,iii^* + 9\,ii^* \\
  & \leq v(D) + 20 + \left(\sum i_n^* + iii^* + ii^*\right) \\ & \leq v(D) + 25. \qedhere
\end{split} 
\]
\endproof

\begin{proposition} \label{prop:square-free}
If $X$ is a square free K3 quartic surface, then $\Phi(X) \leq 55$.
\end{proposition}
\proof
The proof can be copied word by word from \cite[Proposition 5.8]{veniani1}, using both Lemma~\ref{lemma:D-ell} and Lemma~\ref{lemma:D-qe}.
\endproof

\proof[Proof of Theorem \ref{thm:char2}]
Having treated the triangle case (Proposition~\ref{prop:char2-triangle}), the square case (Corollary~\ref{cor:triangle-free+square}) and the square free case (Proposition~\ref{prop:square-free}), the proof is complete.
\endproof

\section{Rams--Schütt's family} \label{sec:char2-family}

The bound of Theorem \ref{thm:char2} is sharp and is reached by all surfaces of the following family, as long as $\lambda \neq 0$:
\[
\XXX: \lambda x_{0} x_{1}^{2} x_{2} + x_{1}^{4} + x_{1} x_{2}^{3} + x_{0}^{3} x_{3} + x_{0} x_{2} x_{3}^{2} = 0
\]

This family was found by Rams and Schütt \cite{char2} and differs by theirs only up to a change of coordinates.

A member $X$ of family $\XXX$, for $\lambda \neq 0$, contains one singular point $P = [0:0:0:1]$ of type $\bA_3$. The point $P$ sits in a configuration $\CCC_1$ lying on the plane $x_0 = 0$. The four lines making up this configuration are cuspidal lines. The remaining $64$ lines -- including, for instance, the line $x_1 = x_3 = 0$ -- are special lines of valency~$19$. The minimal resolution of $X$ is a Shioda supersingular K3 surface of Artin invariant~$2$.

The rest of the section is dedicated to the proof of Theorem~\ref{thm:char2-uniqueness}.

\begin{lemma} \label{lem:C1=>XXX}
    If $\Phi(X) = 68$ and $X$ admits a configuration $\CCC_1$, then $X$ is projectively equivalent to a member of family $\XXX$.
\end{lemma}
\proof
We parametrize the surface in such a way that the configuration $\CCC_1$ sits in the plane $\Pi: x_0 = 0$ and that two lines in $\Pi$ are given by $\ell_1:x_0 = x_1 = 0$ and $\ell_2:x_0 = x_2 = 0$. Let us call the other two lines $\ell_3$ and $\ell_4$. 

Necessarily at least one of the lines in $\Pi$, say, $\ell_1$, has valency greater than $16$. Since $\ell_1$ has degree lower than $3$, it must be a cuspidal line. It follows that the plane $\Pi$ represents a fiber of type $\IV^*$ for $\ell_1$, there are no other singular points on the surface and $v(\ell_1$) is exactly $16$. The same must hold for the other three lines in $\Pi$. Up to coordinate change, we can suppose that one of the lines meeting $\ell_1$ is given by $\ell': x_1 = x_3 = 0$.

Imposing that the lines $\ell_1,\ldots,\ell_4$ are cuspidal, we obtain a quartic which is projectively equivalent to a member of $\XXX$. The following coefficients are different from $0$ and can be normalized to~$1$: 
\[
a_{0220}, a_{1012}, a_{1102}, a_{3001}. 
\]
We first impose that $\ell_1$ and $\ell_2$ are cuspidal setting their polynomial $\varphi$ identically equal to $0$, as in Lemma \ref{lemma:ins=>cusp}, obtaining the following relations:
\begin{gather*}
    a_{1111} = a_{2011} = a_{2020} = a_{2101} = a_{2200} = 0, \\
    a_{0130} = a_{0310} = a_{0220}, \\
    a_{1120} = a_{0130}a_{2002}, \quad a_{1210} = a_{0310} a_{2002}.
\end{gather*}
At this point, $\ell_3$ and $\ell_4$ are given by $x_0 = x_{1}^{2} + x_{1} x_{2} + x_{2}^{2} = 0$
and imposing that they are cuspidal yields the following equation:
\[
     a_{2110} = a_{0220} a_{2002}^{2}.
\]
Changing coordinates
\[
[x_0:x_1:x_2:x_3] \rightsquigarrow [x_0:x_1:a_{2002}x_0 + x_2:(a_{2002}^3+a_{3100})x_1 + x_3],
\]
we recover family $\XXX$.
\endproof

\begin{proposition} \label{prop:triangle=>XXX}
    If $\Phi(X) = 68$ and $X$ contains a triangle, then $X$ is projectively equivalent to a member of family $\XXX$. 
\end{proposition}
\proof
Let $\Pi$ be the plane containing a triangle. By virtue of Proposition~\ref{prop:triangle+sing-pt}, we can suppose that $\Pi$ has no singular points. Let us call $\ell_0,\ldots,\ell_3$ the lines on $\Pi$. At least two of them are special; since each special line induces an automorphism of order $3$ which exchanges the other three, it follows that all four of them induce fibrations with the same singular fiber types and, in particular, the same valency, which must be equal to $19$.

\begin{claim} \label{clm:p-fiber-no-sing-pt}
    A 3-fiber of the lines $\ell_i$, $i = 0,\ldots,3$, cannot contain singular points.
\end{claim}
\proof[Proof of the claim]
By the presence of the automorphism, the residual cubic of the fiber should be as in configuration $\AAA_4$ or $\BBB_4$. We can exclude both of them using formula \eqref{eq:PhiX} and the known bounds of Table \ref{tab:char2-elliptic} (the other three lines must necessarily be elliptic since they admit a fiber of type $\I_n$).
\endproof

Let $\ell'$ be the line in the (ramified) 1-fiber of $\ell_0$ and let $P$ be the point of intersection of $\ell'$ with the residual conic $C$ not on $\ell_0$. Consider the $15$ 3-fibers of $\ell_1$, $\ell_2$ and $\ell_3$ other than $\Pi$. In each of them, there is a line meeting $\ell'$: let us call these lines $m_i$, $i=1,\ldots,15$. By the same argument as before, the lines $m_i$ are special lines of valency~$19$.

We now distinguish the two cases, according to whether $P$ is smooth or singular. Suppose first that $P$ is smooth.

\begin{claim}
    If $P$ is smooth, then $v(\ell') = 16$.
\end{claim}
\proof[Proof of the claim]
Because of the automorphism induced by $\ell_0$, the valency of $\ell'$ has the form $v(\ell') = 1+3\,a$; moreover, $v(\ell') \geq 16$ because $\ell'$ meets $\ell_0$, $m_1,\ldots,m_{15}$, but $v(\ell') \leq 18$ because $\ell'$ is clearly of the first kind. The only possibility is then $v(\ell') = 16$.
\endproof

It follows that $\ell'$ does not belong to a 3-fiber of $m_i$ because otherwise $\ell'$ should have valency $19$; hence, $\ell'$ sits in the 1-fiber of each $m_i$ (and of $\ell_0$), and thus $\ell'$ has exactly $16$ fibers of type $\III$. Necessarily, $\ell'$ is quasi-elliptic, and all other reducible fibers must have an irreducible residual cubic. By an Euler--Poincaré characteristic argument, there can be two cases:
\begin{enumerate}[(i)]
    \item either $\ell'$ has 4 more fibers of type $\III$, whose residual cubics have a cusp which is a singular point of type $\bA_1$;
    \item or $\ell'$ has one fiber of type $\I_0^*$, whose residual cubic is a cusp with a singular point of type $\bD_4$.
\end{enumerate}
\begin{itemize}
    \item In case (i), the line $\ell_0$ has a ramified 1-fiber of type $\I_2$, six 3-fibers without singular point (Claim \ref{clm:p-fiber-no-sing-pt}), and no other 1- or 3-fiber. It follows that $\ell_0$ has four more fibers $F_1,\ldots,F_4$ with irreducible residual cubics, each of them containing one of the four points of type $\bA_1$ (in fact, an irreducible residual cubic can have at most one singular point of the surface). This means that $e(F_i) \geq 2$, so that
    \[
        e(X) \geq 2 + 6\cdot 3 + \sum_{i = 1}^4 e(F_i) \geq 28,
    \]
    which is impossible.
    \item Similarly, in case (ii), the line $\ell_0$ would have one more fiber $F$ containing the point of type $\bD_4$, with $e(F) \geq 6$ (since $F$ contains more than 5 components); we would then obtain
    \[
        e(X) \geq 2 + 6\cdot 3 + e(F) \geq 26, 
    \]
    also impossible.
\end{itemize}

We can therefore suppose that $P$ is a singular point of $X$. Consider now the lines $\ell'_1$, $\ell'_2$, $\ell'_3$ in the 1-fibers of $\ell_1$, $\ell_2$ and $\ell_3$. By the same token, also the line $\ell_i'$ contains a singular point $P_i$, $i=1,2,3$. 
If the four points $P$, $P_1$, $P_2$, $P_3$ are distinct, then $\ell_0$ has a 1-fiber of type $\I_n$, $n\geq 3$, six 3-fibers, and three fibers $F_i$ containing $P_i$ with $e(F_i)\geq 2$, $i=1,2,3$, but this cannot be:
\[
    e(X) \geq 3 + 6\cdot 3 + 3\cdot 2 = 27.
\]
Hence, the points $P_i$ coincide with $P$, and $v(\ell') = 16$. Arguing as before, $\ell'$ has $16$ fibers of type $\III$ and is quasi-elliptic. By Proposition \ref{prop:v>13=>cusp}, $\ell'$ is cuspidal, and so are $\ell_1'$, $\ell_2'$ and $\ell_3'$; necessarily, they must lie in the same plane, forming a configuration $\CCC_1$. We can then conclude applying Lemma \ref{lem:C1=>XXX}.
\endproof

\begin{corollary} \label{cor:comp-red-plane=>XXX}
    If $\Phi(X) = 68$ and $X$ admits a completely reducible plane, then $X$ is projectively equivalent to a member of family $\XXX$.
\end{corollary}
\proof
By Proposition \ref{prop:triangle=>XXX}, we can assume that $X$ is triangle free. By inspection of the proof of Proposition \ref{prop:comp-red-Pi=>leq68}, we see that $X$ must admit a configuration $\CCC_1$, so we can apply Lemma \ref{lem:C1=>XXX}.
\endproof

\proof[Proof of Theorem \ref{thm:char2-uniqueness}]
By Corollary~\ref{cor:comp-red-plane=>XXX}, we can suppose that $X$ does not admit any completely reducible plane. We claim that this assumption leads to a contradiction. 

By virtue of Proposition~\ref{prop:square-free}, we can suppose that $X$ admits a square~$D$ formed, say, by the lines $\ell_1,\ldots,\ell_4$. The square $D$ induces an elliptic fibration~$\pi:X\rightarrow \PP^1$ because $\I_4$ is not a quasi-elliptic fiber. On account of Proposition~\ref{prop:v>13=>leq64}, all lines on $X$ have valency $\leq 13$; hence, by Lemma~\ref{lemma:D-ell} we have
\[
    68 = \Phi(X) \leq 4\cdot(13-2) + 24 = 68.
\]
Since equality holds, we deduce two facts:
\begin{enumerate}[(i)]
    \item $v(\ell_i) = 13$, $i=1,2,3,4$;
    \item all the components of the singular fibers of the fibration~$\pi$ must be lines.
\end{enumerate}
Let $F$ be a general fiber of the fibration $\pi$. Since $F$ is linearly equivalent to $L_1 + L_2 + L_3 + L_4$, we have $F \cdot H = 4$, where $H$ is the hyperplane divisor. It follows from (ii) that all singular fibers are composed by 4 lines. Since the only fiber type with 4 components is $\I_4$, the fibration~$\pi$ has necessarily 6 fibers of type $\I_4$, i.e., 6 squares. Let us put $\ell_i^1 := \ell_i$ and call the other 20 lines $\ell_i^{j}$, $i=1,\ldots,4$, $j=2,\ldots,6$. Arguing as before, we deduce that $v(\ell_i^j) = 13$ for every $i,j$.

Since there are no completely reducible planes, the lines $\ell_i^j$ have no 3-fibers; by Lemma~\ref{lemma:no-p-fibers} and Table~\ref{tab:char2-qe}, they must be quasi-elliptic of degree $3$ or $2$.

\begin{claim}
    The surface $X$ is not smooth.
\end{claim}
\proof[Proof of the claim]
If $X$ were smooth, then $\ell_1$ would induce a fibration with 20 fibers of type $\III$ formed by a line and an irreducible conic, so $v(\ell_1)=20$, which contradicts $v(\ell_1)=13$.
\endproof

Let $E$ be an (irreducible) exceptional divisor coming from the resolution of a singular point $P$. Since $E$ is not a component of a singular fiber of $\pi$, $E$ is a multisection of $\pi$. On the other hand, $E$ can only have positive intersection with one of the $L_i$, $i=1,2,3,4$, since the point $P$ can sit only on one line $\ell_i$, so $E$ is actually a section. Up to index permutation, we can suppose then that $E\cdot L_1^j = 1$ for $j=1,\ldots,6$, i.e., the point $P$ belongs to the lines $\ell_1^j$, $j=1,\ldots,6$ (which have then degree $2$).

\begin{claim}
    The point $P$ is of type $\bA_1$.
\end{claim}
\proof[Proof of the claim]
    If $P$ is not of type $\bA_1$, then there is another exceptional divisor $E'$ coming from $P$, and we can suppose $E\cdot E' = 1$. Arguing as before, $E'$ is also a section and since $P$ is contained in $\ell_1^j$, we have $E'\cdot L_1^j = 1$, too. But $E$ and $E'$ can intersect each $L_1^j$ only in one point, namely the one mapping to $P$ through the resolution, because $\ell_1^j$ has degree $2$ for each $j$. Hence, $E$ and $E'$ have six different points in common, which is impossible since $E\cdot E' = 1$.
\endproof

Finally, let us consider a line in the square $D$ different from $\ell_1$ intersecting~$\ell_1$; up to renaming, we can suppose it is $\ell_2$. Let $\Pi$ be the plane containing both $\ell_1$ and $\ell_2$ and let $F$ be the fiber corresponding to $\Pi$ in the fibration induced by $\ell_2$. The residual conic $C$ in $\Pi$ is irreducible since there are no completely reducible planes. The fiber $F$ is composed of the exceptional divisors coming from $P$, and the strict transforms of $\ell_1$ and $C$, that is to say, three components in total, since $P$ is of type $\bA_1$. On the other hand, $\ell_2$ is a quasi-elliptic line and in characteristic~$2$ there are no quasi-elliptic fibers with three components: contradiction.
\endproof


\printbibliography 

\end{document}